\documentclass[twoside,requo,12pt]{amsart}
\usepackage{amsmath,amsfonts,amssymb}
\usepackage{times}
\usepackage{a4wide}
\usepackage{graphicx}
\usepackage[arrow, matrix, curve, 2cell]{xy}
\usepackage{tikz}
\usepackage{harvard}
\usepackage{hyperref}

\newcounter{mycnt}
\def\themycnt{\thesection.\arabic{mycnt}}
\def\mybenv#1{\refstepcounter{mycnt}%
       \vskip 3pt\noindent{\bf #1~~\themycnt}:~}
\def\myeenv{\hfill\rule{1ex}{1ex}\vskip 3pt}

\makeatletter
\@addtoreset{equation}{section}
\@addtoreset{mycnt}{section}
\makeatother
\def\nn{\nonumber \\}
\def\Id{\text{I\!d}}
\def\calE{\mathcal{E}}

\def\calM{\mathcal{M}}

\def\sfH{\mathsf{H}}
\def\sfR{\mathsf{R}}
\def\sfT{\mathsf{T}}
\def\antip{{\textsf S}}
\def\Bil{\mathsf{Bil}^r_{ass}}
\def\eval{\textsf{ev}}
\def\ceval{\textsf{cev}}
\def\la{\langle}
\def\ra{\rangle}
\def\LR{\mathbf{l}}
\def\RR{\mathbf{r}}
\def\!{\kern -0.15ex}
\def\barbeta{\overline{\beta}}
\def\baralpha{\overline{\alpha}}
\def\Hom{\mathrm{Hom}}
\def\End{\mathrm{End}}
\def\Aut{\mathrm{Aut}}
\def\Der{\mathrm{Der}}
\def\InnDer{\mathrm{InnDer}}
\def\Ext{\mathrm{Ext}}
\def\Ker{\mathrm{Ker}}
\def\Pic{\mathrm{Pic}}

\usetikzlibrary{arrows,backgrounds}
\usetikzlibrary{decorations.pathreplacing}
\usetikzlibrary{decorations.markings}
\pgfdeclarelayer{foreground}
\pgfdeclarelayer{background}
\pgfsetlayers{main,foreground,background}
\newenvironment{pic}[1][]%
{\begin{aligned}\begin{tikzpicture}[#1]}%
{\end{tikzpicture}\end{aligned}}
\tikzstyle{string}=[line width=1.25pt]
\tikzstyle{thickstring}=[line width=2.5pt]
\tikzstyle{dot}=[circle, draw=black, fill=gray, inner sep=.4ex, line width=1.25pt]
\tikzstyle{whitedot}=[circle, draw=black, fill=white, inner sep=.4ex, line width=1.25pt]
\tikzstyle{blackdot}=[circle, draw=black, fill=black, inner sep=.4ex, line width=1.25pt]
\tikzstyle{coupon}=[rectangle,draw=black,minimum size=10pt,line width=1.25pt]
\tikzstyle{cross}=[preaction={draw=white, -, line width=6pt}]
\tikzstyle{mydashed}=[thick,black,dashed] 
\tikzset{arrow/.style={decoration={
    markings,
    mark=at position #1 with \arrow{>}},
    postaction=decorate}
}
\tikzset{reverse arrow/.style={decoration={
    markings,
    mark=at position #1 with \arrow{<}},
    postaction=decorate}
}
\newlength\stateheight
\newlength\minimumstatewidth
\newif\ifhflip\pgfkeys{/tikz/hflip/.is if=hflip}
\makeatletter
\pgfdeclareshape{triangle}
{
    \savedanchor\centerpoint
    {
        \pgf@x=0pt
        \pgf@y=0pt
    }
    \anchor{center}{\centerpoint}
    \anchorborder{\centerpoint}
    \saveddimen\overallwidth
    {
        \pgf@x=3\wd\pgfnodeparttextbox
        \ifdim\pgf@x<\minimumstatewidth
            \pgf@x=\minimumstatewidth
        \fi
    }
    \savedanchor{\upperrightcorner}
    {
        \pgf@x=.5\wd\pgfnodeparttextbox
        \pgf@y=.5\ht\pgfnodeparttextbox
        \advance\pgf@y by -.5\dp\pgfnodeparttextbox
    }
    \anchor{A}
    {
        \pgf@x=-\overallwidth
        \divide\pgf@x by 4
        \pgf@y=0pt
    }
    \anchor{B}
    {
        \pgf@x=\overallwidth
        \divide\pgf@x by 4
        \pgf@y=0pt
    }
    \anchor{text}
    {
        \upperrightcorner
        \pgf@x=-\pgf@x
        \ifhflip
            \pgf@y=-\pgf@y
            \advance\pgf@y by 0.4\stateheight
        \else
            \pgf@y=-\pgf@y
            \advance\pgf@y by -0.4\stateheight
        \fi
    }
    \backgroundpath
    {
        \begin{pgfonlayer}{foreground}
        \pgfsetstrokecolor{black}
        \pgfsetlinewidth{1.25pt}
        \pgfpathmoveto{\pgfpoint{-0.5*\overallwidth}{0}}
        \pgfpathlineto{\pgfpoint{0.5*\overallwidth}{0}}
        \ifhflip
            \pgfpathlineto{\pgfpoint{0}{\stateheight}}
        \else
            \pgfpathlineto{\pgfpoint{0}{-\stateheight}}
        \fi
        \pgfpathclose
        \pgfusepath{stroke}
        \end{pgfonlayer}
    }
}
\makeatother

\begin{document}
\title{Some graphical aspects of Frobenius structures}
%
\author{Bertfried Fauser}
\address{%
	School Of Computer Science\\
	The University of Birmingham\\
	Edgbaston, Birmingham, B25 2TT}
\email{b.fauser@cs.bham.ac.uk}
\thanks{to appear in `Categorical information flow in physics and linguistics',
E. Grefenstette, Ch. Heunen, and M. Sadrzadeh eds.}
%
\subjclass[2000]{Primary
               16W30;		
               Secondary
               18D10;   	
               16T20;		
               16S50;   	
               16H05;  		
               81R50;  		
               }
\keywords{Frobenius structures, graphical calculus, closed structures, yanking}
\date{February 27, 2012}
\begin{abstract}
  We survey some aspects of Frobenius algebras, Frobenius structures
  and their relation to finite Hopf algebras using graphical calculus.
  We focus on the `yanking' moves coming from a closed structure in
  a rigid monoidal category, the topological move, and the `yanking'
  coming from the Frobenius bilinear form and its inverse, used e.g.
  in quantum teleportation. We discus how to interpret the associated
  information flow. Some care is taken to cover non-symmetric Frobenius
  algebras and the Nakayama automorphism. We review graphically the
  Larson-Sweedler-Pareigis theorem showing how integrals of finite
  Hopf algebras allow to construct Frobenius structures. A few pointers
  to further literature are given, with a subjective tendency to
  graphically minded work.
\end{abstract}
\maketitle
\section{Introduction}
\subsection{Scope of this chapter}
  Frobenius algebras surface at many places in mathematics and physics.
Quite recently, using a convenient graphical notation, Frobenius algebras
have been used to investigate foundational issues of quantum theory -- 
references will be given below. Also, as shown elsewhere in this book,
Frobenius algebras emerge in the semantic analysis of natural languages.
The aim of this chapter is to present the basic results about Frobenius
algebras, their relation to finite dimensional Hopf algebras with
special emphasis on using graphical notation. Frobenius structures,
which are related to ring theory, need to be considered too. Frobenius
structures encode such notions as semi simplicity and separability of
rings. We need occasionally to extend the graphical calculus to encode
properties of underlying rings, but will not venture properly into 2-categorical
notions. Moreover, no new results may be found in this chapter, and
far from everything that is known about the subject is covered here. 
However, in passing we will give some pointers to the literature,
which unfortunately is by far to large to be considered completely.
Note that references do not indicate an attribution, but unless
otherwise stated we merely give the source we use.

\subsection{Frobenius' problem}
\index{Frobenius' problem}
In the late nineteenth century Ferdinand August
Frobenius (1849--1917)
\index{Frobenius, August Ferdinand (1849--1917)}
-- and his student Issai Schur-- studied the
representations and characters of the symmetric groups $S_n$. Together
with work by the English mathematicians, notably Alfred Young, this led
to a break through in finite group theory. In early literature,
e.g.~\cite{brauer:nesbitt:1937a,nesbitt:1938a,littlewood:1940a}, the group
algebra $\mathbb{C}[G]$ of a finite group $G$ is synonymously called
`Frobenius algebra'. A finite group has finite order $\vert G\vert$ and one
can form the free $\Bbbk$-vector space $M$ over (the set underlying) $G$. The
group structure then induces a left and right $G$ action on the bimodule $M$
from the algebra structure of $\Bbbk[G]$.
\mybenv{Definition}(regular representations)
\index{representation!regular}
Let $G$ be a finite group, $x_i\in G$, with multiplication $x_i x_j =
\sum_k f^k_{ij} x_k$, with $f^k_{ij}\in\{0,1\}\subset\Bbbk$, and
multiplication table $[f^k_{ij}]$. We associate the following left $\LR$ 
/ right $\RR$ representations on ${}_AA$ / $A_A$ and parastrophic matrix
$\mathbf{p}_{(a)}$ to $A=\Bbbk[G]$.
\vskip1ex

\begin{tabular}[b]{lll}
   $\LR : A \rightarrow \End_\Bbbk({}_A A)$
 & $\LR_{x_i} \cong [f_{i}]_j^k = [f^k_{ij}]$
 & `group matrix'
 \index{group matrix}\index{matrix!group}
\\
   $\RR : A \rightarrow \End_\Bbbk(A_A)$
 & $\RR_{x_j} \cong [f_{j}]_i^k = [f^k_{ij}]$
 & `antistrophic matrix'
 \index{antistrophic matrix}\index{matrix!antistrophic}
\\
 \label{parastrophic}
   $\mathbf{p}_{(a)} : A\otimes A \rightarrow \Bbbk$
 & $\sum_k [f^k_{ij}]a_k = [(P_{(a)})_{ij}]$
 & `parastrophic matrix' ($a_k \in \Bbbk$)
 \index{parastrophic matrix}\index{matrix!parastrophic}
\end{tabular}
\myeenv
It is easy to see that $\LR,\RR$ are, in general reducible, representations
induced by left and right multiplication. The parastrophic matrix is not
a representation, but is related to a linear form on $A$. It contains the
following important information, solving Frobenius' problem of determining
when $\LR$ and $\RR$ are equivalent:
\mybenv{Theorem}\cite{frobenius:1903a}
If there exist $a_k\in \Bbbk$ such that the parastrophic matrix $[P_{(a)}]$ is
invertible, then the left and right regular representations are isomorphic
${}_A A \cong A_A$.
\myeenv
Extending to algebras, we have the
\mybenv{Definition}
An algebra $A$ is called Frobenius iff left and right regular
representations are isomorphic.
\myeenv
\mybenv{Example}
\index{Frobenius algebra!examples}
The reader may check that the commutative polynomial
rings $\Bbbk[X,Y]/$ $\la X^2,Y^2\ra$ and $\Bbbk[X]/\la X^2+1\ra$ are
Frobenius, while $\Bbbk[X,Y]/\la X^2,XY^2,Y^3\ra$ is not Frobenius.
Further examples for Frobenius algebras are the matrix algebras
$A=M_n(\Bbbk)$, where $\Bbbk$ is a division ring. In particular
for $G$ a finite group $A= \mathbb{C}[G]$ is Frobenius.
\myeenv

\subsection{Finite dimensional Hopf algebras}
Studying the topology of group manifolds Heinz Hopf (1894--1971)
introduced in~\cite{hopf:1941a}
\index{Hopf, Heinz (1894--1971)}
the concept of an `Umkehrabbildung',
\index{Umkehrabbildung}
that is a comultiplication. The history of Hopf
algebras is sketched in~\cite{cartier:2007a}. We only note that in the
old days the term `Hopf algebra' is what is now called `Bialgebra'.
Coalgebra is the categorical dual notion to algebra, that is we have a
vector space $C$, and two structure maps $\Delta : C \rightarrow C\otimes C$,
an associative comultiplication, and $\epsilon : C \rightarrow I$ a counit,
fulfilling the axioms obtained by `reversing arrows' in the respective
diagrams for an algebra~\eqref{sliced-unit}. It is convenient to introduce the
Heyneman-Sweedler~\cite{heyneman:sweedler:1969a,heyneman:sweedler:1970a}
\index{Heyneman-Sweedler index notation}
index notation for comultiplications. On an element $c\in C$ one sets
$\Delta(c) = c_{(1)} \otimes c_{(2)} := \sum_i c_{1i} \otimes c_{2i}$.
\mybenv{Definition}\label{definitionHopf}
\index{Hopf algebra!finite, definition of}
  A finite dimensional Hopf algebra $H$ is the sextuple
  $(H, m, \eta, \Delta, \epsilon,$ $\antip)$ where $H$ is a finite dimensional
  vector space, $m,\eta$ are algebra multiplication and unit,
  $\Delta,\epsilon$ are comultiplication and counit, and $\antip$ is the
  antipode, defined as convolutional inverse of the identity
  $m(\antip\otimes \Id)\Delta = \eta\epsilon = m(\Id\otimes \antip)\Delta$,
  fulfilling the compatibility condition: $\Delta(ab)=\Delta(a)\Delta(b)$,
  see \eqref{defHopf}.
\myeenv
The compatibility relation can be read as `the comultiplication is an algebra
homomorphism' (and \emph{vice versa}). A bialgebra is the above structure
without the antipode map. Any graded connected bialgebra is actually a Hopf
algebra. We will see below, that a Frobenius algebra can be described in a
similar way using  a comultiplication. The Frobenius compatibility law is
different, saying that `the comultiplication respects the module structure'
(and \emph{vice versa}). All this will be more obvious when we have the
graphical notation available.

\section{Graphical calculus}

\subsection{History and informal introduction}
\index{graphical calculus}\index{string diagrams}
We work in a (strict) symmetric monoidal category $\mathcal{C}$, with a tensor
as monoidal structure~\cite{majid:1995a,balakov:kirilov:2001a,street:2007a}.
Mainly we are interested in the case of $\mathbf{finVect}_\Bbbk$, finite
dimensional $\Bbbk$-vector spaces, or categories of finite dimensional representations
$\mathcal{R}$, or categories of projective finitely generated left(/right)
modules ${}_{R}\mathcal{M}$ over a (not necessarily commutative unital) ring
$R$. Category theory~\cite{maclane:1971a} comes with the diagrammar of
commutative diagrams (CDs), where objects are represented as vertices and
morphisms as arrows (directed edges) between them. This is one way to define
categories, see~\cite{lambek:scott:1986a}. Graphical calculus was informally
used for a long time, e.g~\cite{brauer:1937a}. Usually its origin is
attributed to Roger Penrose's seminal paper~\cite{penrose:1971a}. The formal
statement that graphical calculus, also called string diagrams, is a sound
transformation of category theory is given
in~\cite{joyal:street:1988a,joyal:street:1991b}. A main thrust for
developing graphical techniques came from low dimensional topology, that is
knot theory~\cite{kauffman:1991a,turaev:1994a,kassel:1995a,ohtsuki:2002a} and
topological quantum field theory, TQFT, e.g.~\cite{atiyah:1989b,kock:2003a}. A
survey and further literature is in~\cite{seelinger:2011a}. Graphical calculus
is in some sense a (Poincar\'e) dual picture to commutative diagrams, where
morphisms are depicted as labelled vertices (depicted also by boxes called
\emph{coupons}) and objects label the edges connecting them. Such a diagram
is called a \emph{tangle},
\index{tangle}\index{graphical calculus!tangle}
\index{graphical calculus!coupon}
it is a representative of an isotopy class of equivalent such
diagrams. Every (unoriented) cycle in a commutative diagram gives rise to a
tangle equation, which establishes a rewriting rule also called a \emph{move}.
\index{graphical calculus!move}
For example the unit law for an algebra $A$ in the monoidal category
$\mathcal{C}$ has a CD with two triangles, and an equivalent description by
two tangle equations
\begin{align}\label{sliced-unit}
\index{graphical calculus!algebra, definition of}
\vcenter{\hsize=0.4\textwidth\xymatrix{
   I\otimes A
   \ar[r]^{\eta\otimes A}\ar[dr]_{\sim}
  &A\otimes A
   \ar[d]^{m}
  &A\otimes I
   \ar[l]_{A\otimes\eta}
   \ar[ld]^{\sim}
  \\
  &
  A
  &
}}
&&
\begin{pic}[scale=0.7]
  \node (in) at (1,2) {};
  \node (eta) [whitedot] at (-1,1) {$\eta$};
  \node (mid) at (1,1) {};
  \node (mul) [blackdot] at (0,0) [label=south east:$m$]{};
  \node (out) at (0,-1) {};
  \draw[string,out=270, in=0] (in) -- (mid.center) to (mul.east);
  \draw[string,out=270, in=180] (eta) to (mul.west);
  \draw[string] (mul) -- (out);
  \node at (1.5,0.1) {$=$};
  \node (1) at (2,2) {};
  \node (2) at (2,-1) {};
  \draw[string] (1.south) -- (2.north);
  \node at (2.5,0.1) {$=$};
  \node (in2) at (3,2) {};
  \node (eta2) [whitedot] at (5,1) {$\eta$};
  \node (mid2) at (3,1) {};
  \node (mul2) [blackdot] at (4,0) [label=south east:$m$] {};
  \node (out2) at (4,-1) {};
  \draw[string,out=270, in=180] (in2) -- (mid2.center) to (mul2.west);
  \draw[string,out=270, in=0] (eta2) to (mul2.east);
  \draw[string] (mul2) -- (out2);
  \draw[mydashed] (-1.25,1.5) -- (5.25,1.5);
  \draw[mydashed] (-1.25,0.5) -- (5.25,0.5);
  \draw[mydashed] (-1.25,-0.65) -- (5.25,-0.65);
  \draw[string,->] (5.5,1.5) to node[xshift=3mm,rotate=90] {time} (5.5,-0.5);
\end{pic}
\end{align}
Here we used $A$ in the CD to denote the identity map $1_A$, and we dropped
the edge label $A$ in the tangles. The tangle was sliced by horizontal dotted
lines so that in each slice only one non-identity operation is performed
(Morse decomposition).
\index{graphical calculus!Morse decomposition}
At this time we also need to make clear, that we read
tangles downward using the \emph{pessimistic arrow of time}
\index{arrow of time!pessimistic}
(Oziewicz, talk at ICCA5, Ixtapa, 1999). Also if (nontrivial) crossings
occur, we use the \emph{left handed} crossings, see~\eqref{tang-example2}.
The reader needs to exercise caution comparing tangle diagrams as some
people are right handed optimists, in that case~\eqref{sliced-unit}
to be read upwards would describe the counit
(relabel: $m\mapsto \Delta$, $\eta\mapsto \epsilon$). We find that
reversing arrows in a CD, and relabelling them, is equivalent to changing
the reading order of the tangle.

Graphical calculus can be interpreted as a `language' built out of basic
\emph{letters}, which form \emph{words} by horizontally (tensoring) or vertically
(composition of morphisms) composing them to form larger tangles. The moves
identify different such words into equivalence classes. Sometimes it is
convenient to introduce special tangles replacing the coupons depicting
them for simplicity and clarity. A selection of graphical entities we are
going to use is given as follows:
\begin{align}\label{tang-example1}
\begin{pic}
  \node (in) at (0,1) {$A$};
  \node[rectangle,draw=black] (mid) at (0,0) {$1_A$};
  \node (out) at (0,-1) {$A$};
  \draw[string] (in.south) to (mid.north);
  \draw[string] (mid.south) to (out.north);
\end{pic}
\,\,&=\hskip-0.5ex
\begin{pic}
  \node (in) at (0,1) {$A$};
  \node (out) at (0,-1) {$A$};
  \draw[string] (in.south) to (out.north);
\end{pic}
\,\,;
\quad
\begin{pic}
  \node (in1) at (0,1) {$A$};
  \node (in2) at (1,1) {$A$};  
  \node[coupon,minimum width=1.3cm] (mid) at (0.5,0) {$m_A$};
  \node (out) at (0.5,-1) {$~A~$};
  \draw[string] (in1.south) to (in1 |- mid.north west);
  \draw[string] (in2.south) to (in2 |- mid.north east);
  \draw[string] (mid.south) to (out.north);
\end{pic}
\,\,=
\begin{pic}
  \node (in1) at (0,1) {$A$};
  \node (in2) at (1,1) {$A$};
  \node[blackdot] (mul) at (0.5,0) {};
  \node (out) at (0.5,-1) {$A$};
  \draw[string,out=270,in=180] (in1) to (mul);
  \draw[string,out=270,in=0] (in2) to (mul);
  \draw[string] (mul) to (out);
\end{pic}
\,\,;
\quad
\begin{pic}
  \node (in1) at (0,1) {$A$};
  \node (in2) at (1,1) {$A$};
  \node[coupon,minimum width=1.25cm] (sig) at (0.5,0) {$\sigma $};
  \node (out1) at (0,-1) {$A$};
  \node (out2) at (1,-1) {$A$};
  \draw[string] (in1) to (in1 |- sig.north west);
  \draw[string] (in2) to (in2 |- sig.north east);
  \draw[string] (out1 |- sig.south west) to (out1);
  \draw[string] (out2 |- sig.south east) to (out2);
\end{pic}
\,\,=
\begin{pic}
  \node (in1) at (0,1) {$A$};
  \node (in2) at (1,1) {$A$};
  \node at (0.8,0.0) {$\sigma $};
  \node (out1) at (0,-1) {$A$};
  \node (out2) at (1,-1) {$A$};
  \draw[string] (in1) to (out2);
  \draw[string] (in2) to (out1);
\end{pic}
\end{align}
We may drop identity morphisms and their coupons. We may use a node to depict
an algebra multiplication (or comultiplication in the inverted reading of the
tangle). We also use the cross for the invertible involutive switch on tensors,
and we drop usually the tensor signs on input and output labels, or even the
labels if they are clear from the context. Creating and deleting elements depict
maps $a : I \rightarrow A$ or $f : A \rightarrow I$, scalars (ring elements)
do not have a graphical counterpart (void diagram) or are depicted by a tangle
with no input and output lines (closed graph). A braiding will be depicted by
keeping over or under information as usually done in knot theory. This reads as
follows:
\begin{align}\label{tang-example2}
\begin{pic}
  \node (out) at (0,-1) {$A$};
  \node[triangle,hflip] (state) at (0,0) {$a$};
  \draw[string] (state) to (out);
\end{pic}
\,\,;
&\quad
\begin{pic}
  \node (in) at (0,1) {$A$};
  \node[triangle] (state) at (0,0) {$f$};
  \draw[string] (in) to (state);
\end{pic}
\,\,;
\begin{pic}
  \node (in1) at (0,1) {$B$};
  \node (in2) at (1,1) {$A$};
  \node[coupon,minimum width=1.25cm] (sig) at (0.5,0) {$R $};
  \node (out1) at (0,-1) {$A$};
  \node (out2) at (1,-1) {$B$};
  \draw[string] (in1) to (in1 |- sig.north west);
  \draw[string] (in2) to (in2 |- sig.north east);
  \draw[string] (out1 |- sig.south west) to (out1);
  \draw[string] (out2 |- sig.south east) to (out2);
\end{pic}
\,\,=
\begin{pic}
  \node (in1) at (0,1) {$B$};
  \node (in2) at (1,1) {$A$};
  \node (out1) at (0,-1) {$A$};
  \node (out2) at (1,-1) {$B$};
  \node at (0.9,0) {$R$};
  \draw[string] (in2) to (out1);
  \draw[string,cross] (in1) to (out2);
\end{pic}
\,\,;
\begin{pic}
  \node (in1) at (0,1) {$B$};
  \node (in2) at (1,1) {$A$};
  \node[coupon,minimum width=1.25cm] (sig) at (0.5,0) {$R^{-1} $};
  \node (out1) at (0,-1) {$A$};
  \node (out2) at (1,-1) {$B$};
  \draw[string] (in1) to ([xshift= 0.15cm] sig.north west);
  \draw[string] (in2) to ([xshift=-0.15cm] sig.north east);
  \draw[string] ([xshift= 0.15cm] sig.south west) to (out1);
  \draw[string] ([xshift=-0.15cm] sig.south east) to (out2);
\end{pic}
\,\,=
\begin{pic}
  \node (in1) at (0,1) {$B$};
  \node (in2) at (1,1) {$A$};
  \node (out1) at (0,-1) {$A$};
  \node (out2) at (1,-1) {$B$};
  \node at (1.1,0) {$R^{-1}$};
  \draw[string] (in1) to (out2);
  \draw[string,cross] (in2) to (out1);
\end{pic}
\end{align}
As in~\eqref{sliced-unit} we sometimes use labeled circles and not triangles to
depict creation and deletion of elements. If we want to distinguish a module
$A$ and its duals $A^*$ (or ${}^*A$) we need either labels or
\emph{oriented tangles}. 
\index{graphical calculus!oriented}\index{tangle!oriented}
We use downward oriented lines for $A$ and upward oriented lines for $A^*$
(and ${}^*A$).

\subsection{Basic rules of graphical calculus} We have not the space to
formally introduce graphical calculus in full detail, so we restrict
ourselves to present the basic facts how to manipulate tangles.

\subsection{Horizontal and vertical composition, sliding} We work in a
\index{tangle!composition of}
rigid symmetric monoidal category $\mathcal{C}$, with a `tensor' bifunctor
$\otimes$ as monoidal structure. We have two types of composition of morphisms.
The (partial) composition of morphisms in the category is called
\emph{vertical composition}, it is depicted as `stacking' (compatible) coupons
of morphisms. We allow rewrites of coupons to be merged at their vertical
boundary.
\index{vertical composition}
\begin{align}
\xymatrix{
   A
   \ar[r]^{f}
   \ar@/_1pc/[rr]_{g\circ f}
  &
   B
   \ar[r]^{g}
  &
   C
}
&&&
\begin{pic}[scale=0.75]
  \node (in) at (0,1.5) {$A$};
  \node[coupon,minimum size=0.6cm] (1) at (0,0.6) {$f$};
  \node at (0.7,0) {$B$};
  \node[coupon,minimum size=0.6cm] (2) at (0,-0.6) {$g$};
  \node (out) at (0,-1.5) {$C$};
  \draw[string] (in) to (1) to (2) to (out);
\end{pic}
=\,\,
\begin{pic}[scale=0.75]
  \node (in) at (0,1.5) {$A$};
  \node[coupon,minimum size=0.6cm] (1) at (0,0) {$g\circ f$};
  \node (out) at (0,-1.5) {$C$};
  \draw[string] (in) to (1) to (out);
\end{pic}
\end{align}
More delicate is the \emph{horizontal composition} of morphisms on tensor
products. We have
\begin{align}
&\vcenter{\hsize=0.4\textwidth{\xymatrix{
  &
   B\otimes C
   \ar[rd]^{B\otimes g}
  &
  \\
  A\otimes C
  \ar[ur]^{f\otimes C}
  \ar[rr]^{f\otimes g}
  \ar[rd]_{A\otimes g}
  &
  &
  B\otimes D
  \\
  &
  A\otimes D
  \ar[ru]_{f\otimes D}
  &
}}}
\qquad\parbox{0.35\textwidth}{with identity maps denoted by $1_A=A$, etc.}
\\
\label{tang-sliding}
&\begin{pic}
  \node (in) at (0,1.5) {$A\otimes C$};
  \node[coupon, minimum width=1cm] (mid) at (0,0) {$f\otimes g$};
  \node (out) at (0,-1.5) {$B\otimes D$};
  \draw[string] (in) to (mid) to (out); 
\end{pic}
=
\begin{pic}
  \node (in1) at (0,1.5) {$A$};
  \node (in2) at (1,1.5) {$C$};
  \node (out1) at (0,-1.5) {$B$};
  \node (out2) at (1,-1.5) {$D$};
  \node[coupon] (ul) at (0,0.5) {$f$};
  \node[coupon] (dr) at (1,-0.5) {$g$};
  \draw[string] (in1) to (ul) to (out1);
  \draw[string] (in2) to (dr) to (out2);
  \draw[mydashed] (-0.3,0.1) rectangle (1.3,0.9);
  \draw[mydashed] (-0.3,-0.9) rectangle (1.3,-0.1);
\end{pic}
\,\,=
\begin{pic}
  \node (in1) at (0,1.5) {$A$};
  \node (in2) at (1,1.5) {$C$};
  \node (out1) at (0,-1.5) {$B$};
  \node (out2) at (1,-1.5) {$D$};
  \node[coupon] (ul) at (0,-0.5) {$f$};
  \node[coupon] (dr) at (1, 0.5) {$g$};
  \draw[string] (in1) to (ul) to (out1);
  \draw[string] (in2) to (dr) to (out2);
  \draw[mydashed] (-0.3,0.1) rectangle (1.3,0.9);
  \draw[mydashed] (-0.3,-0.9) rectangle (1.3,-0.1);
\end{pic}
\,\,=
\begin{pic}
  \node (in1) at (0,1.5) {$A$};
  \node (in2) at (1,1.5) {$C$};
  \node (out1) at (0,-1.5) {$B$};
  \node (out2) at (1,-1.5) {$D$};
  \node[coupon] (ul) at (0,0) {$f$};
  \node[coupon] (dr) at (1,0) {$g$};
  \node at (0.5,0) {$\otimes $};
  \draw[string] (in1) to (ul) to (out1);
  \draw[string] (in2) to (dr) to (out2);
  \draw[mydashed] (-0.3,-0.4) rectangle (1.3,0.4);
\end{pic}
\,\,;
\qquad
\begin{pic}
  \node (in1) at (0,1.5) {$A$};
  \node (in2) at (1,1.5) {$C$};
  \node (out) at (0.5,-1.5) {$B$};
  \node[coupon, minimum width=1.3cm] (tensor) at (0.5,0) {$\otimes $};
  \draw[string] (in1) to (in1 |- tensor.north west);
  \draw[string] (in2) to (in2 |- tensor.north east);
  \draw[string] (tensor) to (out);
\end{pic}
\end{align}
From this we conclude that we are allowed to i) reduce boxes containing
morphisms of the form $f\otimes g$ to two unconnected boxes $f$ and $g$,
and ii) that we are allowed to slide these boxes along each other. Note,
that we drop isomorphisms as shown in the rightmost tangle
in~\eqref{tang-sliding}. This generalizes in an obvious way to $n$
inputs and $m$ outputs of the tangles.

\mybenv{Warning}
If we work in a braided monoidal category, these moves are no longer
available, but need to be altered. Let $R : B\otimes A\rightarrow A\otimes B$
be a braiding, acting on elements as $R(b\otimes a) = a_R\otimes b_R =
a_r\otimes b_r$ using yet a different Heyneman-Sweedler notation for
2-2-morphisms. The ordinary algebra structure of a symmetric monoidal
tensor product on $A\otimes B$ is defined as
$(a\otimes b)(c\otimes d) = ac\otimes bd$ or
$m_{A\otimes B} = (m_A\otimes m_B)(A\otimes \sigma_{A,B}\otimes B)$
where $\sigma$ is the switch map on tensors. Now, let $A\#_R B=A\otimes B$
as a $\Bbbk$-module, and define a new twisted multiplication, the smash
product (recall that $A=1_A$ etc)
\begin{align}
  m_{A\#_R B} &= (m_A\otimes m_B) ( A\otimes R_{A,B}\otimes B)
  &&&
  (a\# b)(c\# d) &= ac_R\# b_Rd\,.
\end{align}
Then the following holds
\index{smash product}
\mybenv{Theorem}\cite[p. 50]{caenepeel:militaru:zhu:2002a} The triple
$(A,B,R)$ is a smash product structure if and only if
\begin{align}\label{smashProduct}
R(b\otimes A) &= A\otimes b
&&&
R(B\otimes a) &= a\otimes B \nonumber\\
R(bd\otimes a) &= a_{Rr}\otimes b_rd_R
&&&
R(b\otimes ac) &= a_Rc_r\otimes b_{Rr}
\end{align}
for all $a,c\in A$ and $b,d\in B$.
\myeenv
The first two identities follow from the unit law, while the second two
relate to (strict) associativity. In this case, we get
$(f\otimes g) = (f\otimes 1)(1\otimes g)$ but the other decomposition is
not direct ($(1\otimes f)(g\otimes 1) = g_R\otimes f_r$). Hence
`sliding' fails to be true and needs modification. More generally
speaking the smash product is related to the question if
$X\cong A\otimes B$ factorizes as an algebra, see \emph{loc. cit.}. There
exists a bijective correspondence between algebra structures on $A\otimes B$
such that the injections $\imath_A, \imath_B$ are algebra maps, and smash product
structures $(A,B,R)$. In what follows we mostly use the graphical
language for the symmetric monoidal case only, hence assume the sliding move
holds.
\myeenv

\subsection{Closed structure, isotopy} A tangle is just a representative
of an equivalence class of diagrams. The equivalence is isotopy of tangles,
which is used in knot theory. We call a tangle an $n$-$m$ tangle, or tangle
with arity $(n,m)$, if it has $n$ input lines (starting at a discrete set of
\index{tangle!arity of}
points on a top horizontal line) and $m$ output lines (ending at a discrete
set of points on a bottom horizontal line), e.g. a multiplication is a 2-1
tangle, a comultiplication is a 1-2 tangle and scalars will be denoted by
0-0 tangles or even by a void tangle. A \emph{critical point} in a tangle
\index{tangle!critical point}
is a point on the plane which has a vertical tangent. We are allowed to
smoothly `bend' lines in a tangle, provided that we keep their topology
and do not introduce or destroy critical points or crossings of lines.
Moves introducing further freedom to modify tangles impose additional
conditions on the underlying category.

In what follows, we need closed structures on the underlying monoidal category.
These come in a left and right version. In the case we have a symmetry $\sigma :
X\otimes Y \rightarrow Y\otimes X$ or a braid, then any two of left duality,
right duality, and braiding defines the third. For an in-depth discussion
see~\cite{kassel:1995a}.
\mybenv{Definition}
\index{monoidal category!rigid}
  A monoidal category $\mathcal{C}$ is \emph{rigid}, if for all $X \in
  \mathcal{C}$ there exist $X^*,{}^*X$ such that the following universal
  morphisms exists:
  \begin{itemize}
  \item Right duality (often denoted also $b_X$, $d_X$):
  \index{duality!right}
  \begin{align}\label{eqRightDual}
     \vcenter{\hsize=0.6\textwidth
       $\begin{array}{ll}
        \eval_X : X^* \otimes X \rightarrow I_X,\quad
        \ceval_X : I_X \rightarrow X \otimes X^*, \\
        \textrm{satisfying}\quad
	  (I_X\otimes \ceval_X)(\eval_X\otimes I_X) = I_X \\
        \phantom{satisfying}\quad
	  (\eval_{X}\otimes I_{X^*})(I_{X^*}\otimes \ceval_{X}) = I_{X^*}
        \end{array}$}
\begin{pic}
  \node (in1) at (0,1) {$X^*$};
  \node (in2) at (1,1) {$X$};
  \node (bot) at (0.5,0) {$\eval_X$};
  \draw[string,arrow=0.55,out=270, in=270] (in2) to (in1);
\end{pic}
;
\begin{pic}
  \node (in1) at (0,-1) {$X$};
  \node (in2) at (1,-1) {$X^* $};
  \node (bot) at (0.5,0) {$\ceval_X$};
  \draw[string,arrow=0.55,out=90, in=90] (in2) to (in1);
\end{pic}
  \end{align}

  \item Left duality (often denoted also $\tilde{b}_X$, $\tilde{d}_X$):
  \index{duality!left}
  \begin{align}\label{eqLeftDual}
     \vcenter{\hsize=0.6\textwidth
       $\begin{array}{ll}
        {}_X\eval : X \otimes X^* \rightarrow I_X,\quad
        {}_X\ceval : I_X \rightarrow X^* \otimes X, \\
        \textrm{satisfying}\quad
          (I_{{^*}X}\otimes {}_{X}\ceval)
          ({}_{X}\eval \otimes I_{{^*}X}) = I_{{^*}X} \\
          \phantom{satisfying}\quad
          ({}_{X}\eval \otimes I_{X})
          (I_{X}\otimes {}_{X}\ceval) = I_{X}
        \end{array}$}
\begin{pic}
  \node (in1) at (0,1) {$X$};
  \node (in2) at (1,1) {${}^*\!X $};
  \node (bot) at (0.5,0) {${}_X\eval$};
  \draw[string,reverse arrow=0.55,out=270, in=270] (in2) to (in1);
\end{pic}
;
\begin{pic}
  \node (in1) at (0,-1) {${}^*\!X $};
  \node (in2) at (1,-1) {$X$};
  \node (bot) at (0.5,0) {${}_X\ceval$};
  \draw[string,reverse arrow=0.55,out=90, in=90] (in2) to (in1);
\end{pic}  \end{align}
  \item Symmetry (braiding) $\sigma_{X,Y} : X\otimes Y \rightarrow Y\otimes X$,
  \index{braid}
    satisfying the braid equation and invertibility
    \begin{align}\label{braid}
     (\sigma\otimes I)(I\otimes \sigma)(\sigma\otimes I) &=
     (\sigma\otimes I)(\sigma\otimes I)(I\otimes \sigma);\qquad
     \sigma \sigma^{-1} = I\otimes I
    \end{align}
  \end{itemize}
  The last equation as tangles represent the Reidemeister 3 move
  \textbf{R3}~\eqref{Reidemeister13} and Reidemeister 2 move
  \textbf{R2}~\eqref{Reidemeister2}.
  \index{Reidemeister moves}\index{graphical calculus!move!Reidemeister}
\myeenv

The conditions on $\eval$, $\ceval$ are depicted as \emph{topological move}
\index{graphical calculus!move!topological}
(or Reidemeister 0 move \textbf{R0}), that is it allows deletion or
introduction of two compatible extrema. This move is also ambiguously
addressed as `yanking',
\index{yanking}
but we will see below that Frobenius bilinear forms also allow a `yanking' of
lines, so we reject this term. We have introduced oriented lines, to depict
objects $X$, $X^*$ and ${}^*X$. If the distinction between left and right
dual is vital, we need to apply labels. The topological move $\textbf{R0}_r$
for the right duality (the left dual tangles are obtained by inverting
orientation) depicts the conditions in~\eqref{eqRightDual}
(and with inverted orientation that of~\eqref{eqLeftDual}).
\begin{align}\label{R0}
  \textbf{R0}_r&:\qquad
\begin{pic}
  \node (in) at (0,1) {};
  \node (1) at (0,0) {};
  \node (2) at (1,0) {};
  \node (3) at (2,0) {};
  \node (out) at (2,-1) {};
  \draw[string] (in) to (1.center);
  \draw[string] (3.center) to (out);
  \draw[string,reverse arrow=0.55,out=270,in=270] (1.center) to (2.center);
  \draw[string,reverse arrow=0.55,out= 90,in= 90] (2.center) to (3.center);
\end{pic}
\,\,=
\quad
\begin{pic}
  \node (in) at (0,1) {};
  \node (out) at (0,-1) {};
  \draw[string,reverse arrow=0.5] (in) to (out);
\end{pic}
\,\,;
\qquad
\begin{pic}
  \node (in) at (0,-1) {};
  \node (1) at (0,0) {};
  \node (2) at (1,0) {};
  \node (3) at (2,0) {};
  \node (out) at (2,1) {};
  \draw[string] (in) to (1.center);
  \draw[string] (3.center) to (out);
  \draw[string,reverse arrow=0.55,out= 90,in= 90] (1.center) to (2.center);
  \draw[string,reverse arrow=0.55,out=270,in=270] (2.center) to (3.center);
\end{pic}
\,\,=
\quad
\begin{pic}
  \node (in) at (0,1) {};
  \node (out) at (0,-1) {};
  \draw[string,arrow=0.55] (in) to (out);
\end{pic}\,.
\end{align}
The Reidemeister 2 move \textbf{R2} depends on what is encoded by the
lines in a tangle, see section~\ref{thick-tangles}. If the lines are
assumed to be one dimensional `strings' (\emph{sic} the name of the
calculus), then straightening a loop introducing a twist $\theta$ in
a string does not matter. In this case the Reidemeister 2 move is
given by (graphically as lhs of \eqref{Reidemeister2})
\begin{align}
 \theta := (I\otimes \eval)(\sigma\otimes I)(I\otimes \ceval) & = I
\end{align}
If we assume that lines in a tangle are more complicated objects,
e.g. ribbons or cylinders (in TQFT, see sec.~\ref{thick-tangles}),
then we need to keep track of the twists. In this case
$\theta \not= I$, and the Reidemeister 2 move needs another loop
$\theta^{-1}$ with an inverse braid to compensate.
\index{twist morphism}

To summarize, we have the following isotopy moves on tangles
\begin{align}\label{Reidemeister13}
&
\begin{pic}
  \node (in1) at (-0.5,1) {};
  \node (in2) at (0.5,1) {};
  \node (out1) at (-0.5,-1) {};
  \node (out2) at (0.5,-1) {};
  \draw[string] (in1) .. controls (0.5,0) and (0.5,0) .. (out1);
  \draw[string,cross] (in2) .. controls (-0.5,0) and (-0.5,0) .. (out2);  
  \node at (0.75,0) {$=$};
  \node (in1) at (1.25,1) {};
  \node (in2) at (1.75,1) {};
  \node (out1) at (1.25,-1) {};
  \node (out2) at (1.75,-1) {};
  \draw[string] (in1) to (out1);
  \draw[string,cross] (in2) to (out2);
  \node[below] at (1.5,-1) {Reidemeister 1};  
  \node at (2.25,0) {$=$};
  \node (in1) at (2.5,1) {};
  \node (in2) at (3.5,1) {};
  \node (out1) at (2.5,-1) {};
  \node (out2) at (3.5,-1) {};
  \draw[string] (in2) .. controls (2.5,0) and (2.5,0) .. (out2);
  \draw[string,cross] (in1) .. controls (3.5,0) and (3.5,0) .. (out1);
\end{pic}
\,\,;
\quad
\begin{pic}
  \node (in1) at (0,1) {};
  \node (in2) at (0.5,1) {};
  \node (in3) at (1,1) {};
  \node (out1) at (0,-1) {};
  \node (out2) at (0.5,-1) {};
  \node (out3) at (1,-1) {};
  \draw[string] (in3) .. controls (1,-0.5) and (0,0) .. (out1);
  \draw[string,cross] (in2) .. controls (0,0.0) and (0,0.0) .. (out2);
  \draw[string,cross] (in1) .. controls (0,0.5) and (1,-0.5) .. (out3);
  \node at (1.25,0) {$=$};
  \node[below] at (1.25,-1) {Reidemeister 3};  
  \node (in4) at (1.5,1) {};
  \node (in5) at (2,1) {};
  \node (in6) at (2.5,1) {};
  \node (out4) at (1.5,-1) {};
  \node (out5) at (2,-1) {};
  \node (out6) at (2.5,-1) {};
  \draw[string] (in6) .. controls (2.5,0.5) and (1.5,0.5) .. (out4);
  \draw[string,cross] (in5) .. controls (2.5,0.0) and (2.5,0.0) .. (out5);
  \draw[string,cross] (in4) .. controls (1.5,0.5) and (2.5,-0.5) .. (out6);
\end{pic}
\end{align}
where the braiding is trivial for the switch map $\sigma$. The two
alternative Reidemeister 2 moves (with and without a twist/braiding) read
\begin{align}\label{Reidemeister2}
&
\begin{pic}
  \node (in) at (0,1) {};
  \node[coupon] (theta) at (0,0) {$\theta $};
  \node (out) at (0,-1) {};
  \draw[string] (in) to (theta) to (out);
  \node at (0.5,0) {$:=$};
  \node (in1) at (1,1) {};
  \node (midn) at (1.5, 0.25) {};
  \node (mids) at (1.5,-0.25) {};
  \node (out1) at (1,-1) {};
  \draw[string,out=180, in=90] (midn.center) to (out1);
  \draw[string,out=0,in=0] (mids.center) .. controls (1.8,-0.25) and (1.8,0.25) .. (midn.center);
  \draw[string,cross,out=270, in=180] (in1) to (mids.center);
  \node[below] at (1.25,-1) {Reidemeister 2};
  \node at (2,0) {$=$};
  \node (in2) at (2.5,1) {};
  \node (out2) at (2.5,-1) {};
  \draw[string] (in2) to (out2);
\end{pic}
\,\,;
\quad
\begin{pic}
  \node (in) at (0,1) {};
  \node (out at (0,-1) {};
     \node (m1) at (0.5, 0.75) {};
     \node (m2) at (0.5, 0.25) {};
     \node (m3) at (0.5,-0.25) {};
     \node (m4) at (0.5,-0.75) {};
  \draw[string,out=0,in=0] (m1.center) .. controls (0.8, 0.75) and (0.8, 0.25) .. (m2.center);    
  \draw[string,out=0,in=0] (m1.center) .. controls (0.0, 0.75) and (0.0,-0.75) .. (m4.center);
  \draw[string,out=0,in=0] (m3.center) .. controls (0.8,-0.25) and (0.8,-0.75) .. (m4.center);   
  \draw[string,cross,out=270, in=180] (in) to (m2.center);
  \draw[string,cross,out=180, in=90] (m3.center) to (out);
  \node[below] at (0.75,-1) {Reidemeister 2'};
  \node at (1,0) {$=$};
  \node (in1) at (1.5,1) {};
  \node (out1) at (1.5,-1) {};
  \draw[string] (in1) to (out1);
\end{pic}
\end{align}
In the symmetric monoidal case the twist morphism $\theta_A$ can be seen
roughly as the composition $A\rightarrow {}^*A \rightarrow ({}^*A)^*$. If
${}^*A\cong A^*$ then $\theta_A$ is the canonical identification
$A\cong (A^*)^*$, explaining why loops can be undone. The modified
Reidemeister 2' move equals identity using the two identifications $A\cong
(A^*)^*$ and $A\cong {}^*({}^*A)$ explaining why two loops are necessary,
and the braiding of course.

Moreover, we can move lines above and below extrema (evaluation and
coevaluations). As an example look at
\begin{align}
&
\begin{pic}[scale=0.7]
  \node (in) at (0,1.5) {};
  \node (out1) at (0,-1.5) {};
  \node (out2) at (1.5,-1.5) {};
  \node (out3) at (2,-1.5) {};
  \draw[string,reverse arrow=0.55] (out1) .. controls (0,1) and (2,2) .. (out3);
  \draw[string,cross,arrow=0.55] (in) .. controls (0,0.5) and (1.5,0) .. (out2);
\end{pic}
=
\begin{pic}[scale=0.7]
  \node (in) at (0,1.5) {};
  \node (out1) at (-0.5,-1.5) {};
  \node (out2) at (1.75,-1.5) {};
  \node (out3) at (2.5,-1.5) {};
  \node (m1) at (-0.5,0) {};
  \node (m2) at (0.5,0) {};
  \node (m3) at (1.5,0) {};
  \node (m4) at (2.5,0) {};
  \draw[string,reverse arrow=0.5] (m3.center) .. controls ( 1.5,1.25) and (2.5,1.25) .. (m4.center);
  \draw[string,reverse arrow=0.55] (m2.center) .. controls ( 0.5,-0.5) and (1.5,-0.5) .. (m3.center);
  \draw[string,reverse arrow=0.55] (m1.center) .. controls (-0.5,0.5) and (0.5,0.5) .. (m2.center);
  \draw[string,arrow=0.55] (m1.center) to (out1);
  \draw[string] (m4.center) to (out3);
  \draw[string,cross, arrow=0.5,out=270,in=90] (in) .. controls (0,0.5) and (1.75,1.5) .. (out2);
\end{pic}
=
\begin{pic}[scale=0.7]
  \node (in) at (2,1.5) {};
  \node (out1) at (0,-1.5) {};
  \node (out2) at (1.5,-1.5) {};
  \node (out3) at (2,-1.5) {};
  \draw[string,reverse arrow=0.5] (out1) .. controls (0,2) and (2,1) .. (out3);
  \draw[string,cross,arrow=0.75] (in) .. controls (2,1) and (1.5,0) .. (out2);
\end{pic}
;\,\,\textrm{with}
\begin{pic}[scale=0.7]
  \node (in1) at (0,1.5) {};
  \node (in2) at (1,1.5) {};
  \node (out1) at (0,-1.5) {};
  \node (out2) at (1,-1.5) {};
  \draw[string,reverse arrow=0.25] (in1) to (out2);
  \draw[string,cross,arrow=0.75] (in2) to (out1);
\end{pic}
:=
\begin{pic}[scale=0.7]
  \node (in1) at (0,1.5) {};
  \node (in2) at (0.5,1.5) {};
  \node (mid) at (1,0) {};
  \node (out1) at (1.5,-1.5) {};
  \node (out2) at (2,-1.5) {};
  \draw[string,reverse arrow=0.8] (in1) .. controls (0,-0.75) and (0.75,-0.75) .. (mid.center);
  \draw[string,reverse arrow=0.3] (mid.center) .. controls (1.25,0.75) and (2,0.75) .. (out2);
  \draw[string,cross,arrow=0.75] (in2) to (out1);
\end{pic}
\end{align}
The other cases are similar.

To relate right and left duality, we need the braiding. If
$\sigma^2=I\otimes I$, is the ordinary switch, we get a symmetric monoidal
category and can drop the left/right distinction, as they are related by
the identity morphism, see left tangle equation in~\eqref{leftFromRight}.
If $\sigma$ is a proper braid~\eqref{braid}, then we can derive the left
duality from the right duality. This is shown in the following equations
for the `cup' tangles $\eval$
\begin{align}\label{leftFromRight}
\begin{pic}[scale=0.7,baseline=0.75cm]
  \node (in1) at (0,1) {};
  \node (in2) at (1,1) {};
  \node (m1) at (0,0) {};
  \node (m2) at (1,0) {};
  \draw[string] (in1) to (m1.center);
  \draw[string] (in2) to (m2.center);
  \draw[string, reverse arrow=0.5] (m2.center) .. controls (1,-0.5) and (0,-0.5) .. (m1.center);
  \node at (1.5,0) {$=$};
  \node[below] at (1.5,-0.55) {if symmetric};
  \node (in3) at (2,1) {};
  \node (in4) at (3,1) {};
  \node (m3) at (2,0) {};
  \node (m4) at (3,0) {};
  \draw[string] (in3) .. controls (2,0.5) and (3,0.5) .. (m4.center);
  \draw[string] (in4) .. controls (3,0.5) and (2,0.5) .. (m3.center);
  \draw[string, arrow=0.55] (m4.center) .. controls (3,-0.5) and (2,-0.5) .. (m3.center);
\end{pic}
;\,\,
\begin{pic}[scale=0.7]
  \node (in1) at (0,1) {};
  \node (in2) at (1,1) {};
  \node (m1) at (0,0) {};
  \node (m2) at (1,0) {};
  \draw[string] (in1) to (m1.center);
  \draw[string] (in2) to (m2.center);
  \draw[string, reverse arrow=0.5] (m2.center) .. controls (1,-0.5) and (0,-0.5) .. (m1.center);
\end{pic}
=
\begin{pic}[scale=0.6]
  \node (in1) at (1,1.5) {};
  \node (in2) at (2,1.5) {};
  \node (l1) at (0,1) {};
  \node (l2) at (0,0.25) {};
  \node (l3) at (0,-0.25) {};
  \node (l4) at (0,-1) {};
  \node (m1) at (1,1) {};
  \node (m4) at (1,-1) {};
  \node (r4) at (2,-1) {};
  \draw[string,->] (in1) to (m1.center);
  \draw[string,reverse arrow=0.75] (in2) to (r4.center);
  \draw[string,arrow=0.5] (l1.center) .. controls (1,0) and (1,0) .. (l4.center);
  \draw[string,cross,out=270,in=0] (m1.center) to (l2.center);
  \draw[string,cross,out=0,in=90] (l3.center) to (m4.center);
  \draw[string,reverse arrow=0.5] (l1.center) .. controls (-0.5,1.5) and (-0.75,0.25) .. (l2.center);
  \draw[string,reverse arrow=0.5] (l3.center) .. controls (-0.75,-0.25) and (-0.5,-1.5) .. (l4.center);
  \draw[string,arrow=0.55,in=270,out=270] (m4.center) to (r4.center);
\end{pic}
=
\begin{pic}[scale=0.6]
  \node (i1) at (1,1.5) {};
  \node (i2) at (2,1.5) {};
  \node (l1) at (0,1) {};
  \node (l2) at (0,0.25) {};
  \node (l3) at (0,-0.25) {};
  \node (l4) at (0,-1) {};
  \node (m1) at (1,1) {};
  \node (m2) at (1,0.25) {};
  \node (m3) at (1,-0.25) {};
  \node (m4) at (1,-1) {};
  \node (r1) at (2,1) {};
  \node (r4) at (2,-1) {};
  \node (rr1) at (2.5,-1) {};
  \node (rr2) at (3,-1) {};
   \draw[string,out=0,in=90] (m2.center) to (rr2.center);
   \draw[string,out=0,in=90] (m3.center) to (rr1.center);
   \draw[string,reverse arrow=0.5] (r4.center) .. controls (2,-1.5) and (1,-1.5) .. (m4.center); 
   \draw[string,arrow=0.5] (rr2.center) .. controls (3,-1.5) and (2.5,-1.5) .. (rr1.center);
   \draw[string,cross,arrow=0.75] (r4.center) to (i2);
   \draw[string,out=0,in=180] (l1.center) to (m2.center);
   \draw[string,cross,out=270,in=0] (m1.center) to (l2.center);
   \draw[string,out=180,in=0] (m3.center) to (l4.center);
   \draw[string,cross,out=0,in=90] (l3.center) to (m4.center);
   \draw[string,->] (i1) to (m1.center);
   \draw[string,reverse arrow=0.5] (l1.center) .. controls (-0.5,1) and (-0.5,0.25) .. (l2.center);
   \draw[string,reverse arrow=0.5] (l3.center) .. controls (-0.5,-0.25) and (-0.5,-1) .. (l4.center);
\end{pic}
=\,\,
\begin{pic}[scale=0.7]
  \node (in1) at (0,1) {};
  \node (in2) at (1,1) {};
  \node[coupon] (theta) at (0,0.25) {$\theta $};
  \node (l2) at (0,-0.75) {};
  \node (r1) at (1,0) {};
  \node (r2) at (1,-0.75) {};
  \draw[string] (in1) to (theta);
  \draw[string] (in2) to (r1.center);
  \draw[string,out=270,in=90] (theta) to (r2.center);
  \draw[string,cross,out=270,in=90] (r1.center) to (l2.center);
  \draw[string,reverse arrow=0.5,out=270,in=270] (l2.center) to (r2.center);  
\end{pic}
\end{align}
The `cap' tangles $\ceval$ are related in a similar fashion. There are some
subtleties in the interplay between dualities and twists, which we do not
contemplate here further,
see~\cite{joyal:street:1991c,turaev:1994a,kassel:1995a,street:2007a}.

In a symmetric monoidal category we can assume \textbf{R0, R1, R2, R3} with a
trivial twist. As we also get ${}^*A\cong A^*$ we can drop orientation too.
This is called `ambient isotopy'. If we work in a ribbon category, we replace
\textbf{R2} by its modified version \textbf{R2}' and keep track of
twist and braid morphisms, this is called `regular isotopy'; for
terminology see~\cite{kauffman:1991a}. In what follows we will for
simplicity work mainly in the symmetric monoidal setting.

\subsection{Tangles not depicting strings}\label{thick-tangles}
\index{tangle!thickened}
`String diagrams', related to symmetric monoidal categories, take their name from
picturing them literally as strings, assumed to have zero radial extension.
Such strings cannot be twisted (or twisting them is irrelevant). However,
there are mathematical structures which are sensitive to twisting. They may
be depicted e.g. by ribbons, and lead accordingly to `ribbon categories'. This
is a reason to speak instead about `tangle diagrams'. However, there are still
more thickenings of tangles, such as cylinders in topological quantum field
theory (TQFT). If we still use `strings' to depict them, we are forced to
adjust the allowed moves, e.g. change the Reidemeister 2 move. This leads to
different notions of isotopy, see~\cite{kauffman:1991a}. Here we just depict
ribbons and the relevant tangles for TQFT for reference, but in the sequel
we will just use strings which reduces the artwork considerably and does not
loose information if care is taken about using only allowed rewriting rules.
\begin{align}
\begin{pic}
  \node (in) at (0,1.1) {};
  \node (out) at (0,-1.1) {};
  \node (m1) at (0.6,0.5) {};
  \node (m2) at (0.6,-0.5) {};
  \draw[line width=1.25pt,draw=black,double=white,double distance=5pt,out=180,in=90] (m1.center) to (out);  
  \draw[line width=1.25pt,draw=black,double=white,double distance=5pt,out=270,in=180] (in) to (m2.center);
  \draw[line width=1.25pt,draw=black,double=white,double distance=5pt] (m1.center) .. controls (1.25,0.5) and (1.25,-0.5) .. (m2.center);
\end{pic}
\,\,=
\quad
\begin{pic}
  \node (in) at (0,1.1) {};
  \node (il) at (-3pt,1.1) {};
  \node (ir) at ( 3pt,1.1) {};
  \node (out) at (0,-1.1) {};
  \node (ol) at (-3pt,-1.1) {};
  \node (or) at ( 3pt,-1.1) {};
  \node (m1) at (0,0.5) {};
  \node (m2) at (0,-0.5) {};
  \draw[mydashed] (in) to (out);
  \draw[string,out=270,in=90+30] (il) to (m1.center);
  \draw[string,out=270,in=90-30] (ir) to (m1.center);
    \filldraw[fill=gray,string,out=270-30,in=90+30] (m1.center) to (m2.center);
    \shadedraw[shading=axis,shading angle=90,fill=gray,string,out=270+30,in=90-30] (m1.center) to (m2.center);
  \draw[string,out=270-30,in=90] (m2.center) to (ol);
  \draw[string,out=270+30,in=90] (m2.center) to (or);
  \draw[string,->] (-0.4,0.1) arc (215:330:0.5cm); 
\end{pic}
\,\,;\quad
\begin{pic}[scale=0.7]
  \draw[line width=1.25pt] (0,1) ellipse (0.5cm and 0.25cm);
  \draw[line width=1.25pt] (-0.5,1) .. controls (-0.5,-0.5) and (0.5,-0.5) .. (0.5,1);
\end{pic}
\,\,;
\quad
\begin{pic}[scale=0.7]
  \begin{scope}
    \clip (0.2,-1) rectangle (1.3,-0.5);
    \draw[mydashed] (0.75,-1) ellipse (0.5cm and 0.25cm);
  \end{scope}\begin{scope}
    \clip (0.2,-1) rectangle (1.3,-1.5);
    \draw[line width=1.25pt] (0.75,-1) ellipse (0.5cm and 0.25cm);
  \end{scope}
  \draw[line width=1.25pt] (0,1) ellipse (0.5cm and 0.25cm);
  \draw[line width=1.25pt] (1.5,1) ellipse (0.5cm and 0.25cm);
  \draw[line width=1.25pt] (0.5,1) .. controls (0.5,0) and (1,0) .. (1,1);
  \draw[line width=1.25pt] (-0.5,1) .. controls (-0.5,0) and (0.25,0) .. (0.25,-1);
  \draw[line width=1.25pt] (2,1) .. controls (2,0) and (1.25,0) .. (1.25,-1);
\end{pic}
\,\,;
\quad
\begin{pic}[scale=0.7]
  \draw[line width=1.25pt] (0,1) ellipse (0.5cm and 0.25cm);
  \begin{scope}
    \clip(-0.6,-1) rectangle (0.6,-0.5);
    \draw[mydashed] (0,-1) ellipse (0.5cm and 0.25cm);
  \end{scope}
  \draw[string] (-0.5,1) to (-0.5,-1);
  \draw[string] (0.5,1) to (0.5,-1);
  \draw[string] (0,0.75) .. controls (0,0.5) and (0.2,0.5) .. (0.5,0.25);
  \draw[mydashed] (0.5,0.25) .. controls (0,0.1) and (0,-0.1) .. (-0.5,-.25);
  \draw[string] (0,-1.25) .. controls (0,-0.5) and (-0.2,-0.5) .. (-0.5,-.25);
  \clip (-0.6,-1.5) rectangle (0.6,-1);
  \draw[line width=1.25pt] (0,-1) ellipse (0.5cm and 0.25cm);
\end{pic}
\end{align}
\index{tangle!ribbon}
\index{topological quantum field theory (TQFT)!disk, trinion, cylinder}
The left equation shows how a loop, if straightened, produces a $2\pi$ twist
in the ribbon. The three right most diagrams depict two dimensional surfaces,
cobordisms, which connect the (oriented, one dimensional) circles at the top
with (oriented) circles at the bottom. They are called `disk' (the unit, or
counit if inverted), `trinion' (the Frobenius algebra product, or coproduct
if inverted), and the `cylinder' (identity map). We have depicted a twisted
line on the cylinder, showing that one can also have a twist (Dehn twist) on
cylinders. A TQFT based on these diagrams allows all two dimensional
Riemannian surfaces to be constructed and characterized, see~\cite{lawrence:1996a}.

\medskip
\section{Frobenius and Hopf algebras}

  In this section we recall some facts about, and most importantly some
characterizations of, Frobenius algebras. We want to emphasize especially
the difference between the multiplications in endormorphisms rings
$\End_R(M)$ over an $R$-module $M$ and the left (right) action
induced by an algebra multiplications $m_A : A \otimes A \rightarrow A$,
where we look at the right (left) factor $A$ as a left $A$-module ${}_AA$
(right $A$-module $A_A$). Our main sources
are~\cite{yamagata:1996a,kadison:1999a,caenepeel:militaru:zhu:2002a,murray:2005a,lorenz:2011a,lorenz:tokoly:2010a}.
General texts on Hopf algebras and modules
are~\cite{sweedler:1969a,abe:1980a,kasch:1982a,caenepeel:1998a,street:2007a}.

  Usually we work over a field, but we will recall here some more general
notions where we also allow more general base rings $R$, examples being
a residual field, a ring extension or even a noncommutative ring.

\subsection{Actions, coactions, representations and two multiplications}
  To understand the similarities and differences between the Frobenius and
endormorphism structures, we need to look briefly at algebra representations.
\index{algebra!representation}
We do that superficially only to the extent which is necessary for our purpose.

  Let ${}_RM_S$ be a (p.f.) $R,S$-bimodule and $\calE(M) = \End(M)
\cong M \otimes_S M^*$ be the endormophism ring over $M$. Let $\{x_i\}_{i=1}^n$
be a set of generators for $M$ and $\{f_i\}_{i=1}^n$ be a dual basis for
$M^*$, i.e. $\eval(f_i \otimes x_j)=f_i(x_j)=\delta_{ij}$.
Given a nondegenerate bilinear form $\beta : M\otimes M\rightarrow R$ with
inverse $\barbeta=\sum x_i \otimes y_i$ (see~\eqref{betaInv}). The bilinear
form provides us with another set of generators $\{y_i\}_{i=1}^n$ for
$M$, such that $\sum\beta(m,x_i)y_i = m$ for all $m\in M$.
For simplicity we denote $\beta^*=\barbeta$, as it is distinguished by
its type signature, see~\eqref{tangle-FrobIso}. The next tangles describe
a left $A$ action on ${}_RM_S$, and how left-duality allows to define therefrom
a right $A$ action on ${}_SM^*_R$. For the moment we use 2-tangles, where the 
area depicts the ring in question, for more details on such tangles see
\index{tangle!2-categorical}
e.g.~\cite{lauda:2006a,khovanov:2010a}. The rightmost tangle is a coaction
for which similar results hold by tangle symmetry.
\begin{align}\label{leftrightaction}
\begin{pic}
  \node (in1) at (0,0.75)  [label=above:$A$] {};
  \node (in2) at (0.75,0.75)  [label=above:${}_R M_S$] {};
  \node (out) at (0.75,-0.75) [label=below:${}_R M_S$] {};
  \node (leftact) at (0.75,0) {};
  \draw[thickstring,arrow=0.8] (in2.center) to node[auto]{$S$} (out.center);
  \draw[string] (in1.center) -- node[auto,swap]{$R$} node[auto]{$R$} (in1 |- leftact) -- (leftact.center);
\end{pic}
\,;\,
\begin{pic}
  \node (in1) at (0.75,0.75)  [label=above:$A$] {};
  \node (in2) at (0,0.75)  [label=above:${}_S M^*_R$] {};
  \node (out) at (0,-0.75) [label=below:${}_S M^*_R$] {};
  \node (leftact) at (0,0) {};
  \draw[thickstring,reverse arrow=0.3] (in2.center) to node[auto,swap]{$S$} (out.center);
  \draw[string] (in1.center) -- node[auto,swap]{$R$} node[auto]{$R$} (in1 |- leftact) -- (leftact.center);
\end{pic}
\hskip-1ex:=\hskip-1ex
\begin{pic}
  \node (in1) at (0,0.75)  [label=above:${}_S M^*_R$] {};
  \node (in2) at (0.75,0.75)  [label=above:$A$] {};
  \node (out) at (2.5,-0.75) [label=below:${}_S M^*_R$] {};
  \node (leftact) at (1.5,0) {};
  \node (m1) at (0,-0.25) {};
  \node (m2) at (1.5,-0.25) {};
  \node (m3) at (1.5,0.25) {};
  \node (m4) at (2.5,0.25) {};  
  \draw[thickstring,reverse arrow=0.3] (in1.center) to node[auto,swap]{$S$} (m1.center);
  \draw[thickstring,out=270,in=270] (m1.center) to (m2.center);
  \draw[thickstring,out=90,in=90] (m3.center) to (m4.center);
  \draw[thickstring] (m3.center) to (m2.center);
  \draw[thickstring,reverse arrow=0.8] (m4.center) to node[auto]{$R$} node[auto,swap]{$S$} (out.center);
  \draw[string] (in2.center) -- node[auto,swap]{$R$} node[auto]{$R$} (in2 |- leftact) -- (leftact.center);
\end{pic}
\,;\,
\begin{pic}
  \node (in1) at (0,-0.75)  [label=below:$A$] {};
  \node (in2) at (0.75,-0.75)  [label=below:${}_R M_S$] {};
  \node (out) at (0.75,0.75) [label=above:${}_R M_S$] {};
  \node (leftact) at (0.75,0) {};
  \draw[thickstring,reverse arrow=0.3] (in2.center) to node[auto,swap]{$S$} (out.center);
  \draw[string] (in1.center) -- node[auto,swap]{$R$} node[auto]{$R$} (in1 |- leftact) -- (leftact.center);
\end{pic}
\end{align}
  The Frobenius property induces an isomorphism between left modules ${}_AM$ and right
modules $M_A$, which does not follow from duality alone. We can use the bilinear from
to define a right action on $M$ from the right action on $M^*$ as follows
\index{tangle!Frobenius isomorphism}
\begin{align}\label{tangle-FrobIso}
\begin{pic}
  \node (in1) at (0,0.75) {};
  \node (in2) at (0.75,0.75) {};
  \node (out) at (0.75,-0.75) {};
  \node (leftact) at (0.75,0) {};
  \draw[thickstring,arrow=0.8] (in2.center) to (out.center);
  \draw[string] (in1.center) -- (in1 |- leftact) -- (leftact.center);
\end{pic}
\stackrel{\text{Frobenius}}{\cong}
\begin{pic}
  \node (in1) at (0.75,0.75) {};
  \node (in2) at (0,0.75) {};
  \node (out) at (0,-0.75) {};
  \node (leftact) at (0,0) {};
  \draw[thickstring,arrow=0.8] (in2.center) to (out.center);
  \draw[string] (in1.center) -- (in1 |- leftact) -- (leftact.center);
\end{pic}
:=
\begin{pic}
  \node (in1) at (0,0.75) {};
  \node (in2) at (1.625,0.75) {};
  \node[coupon,minimum width=1cm] (b) at (0.375,-0.2) {$\beta $}; 
  \node[coupon,minimum width=1cm] (bb) at (1.625,-0.25) {$\barbeta $};
  \node (m1) at (0.75,0.5) {};
  \node (m2) at (1.25,0.5) {};
  \node (m3) at (1.625,0.3) {};
  \node (m4) at (1.95,0.25) {};
  \node (m5) at (2.5,0.25) {};
  \node (o) at (2.5,-0.75) {};
  \draw[thickstring,arrow=0.6,out=90,in=90] (m4.center) to (m5.center);
  \draw[thickstring] (m4.center) to (m4 |- bb.north);
  \draw[thickstring] (m5.center) to (o.center);
  \draw[thickstring,arrow=0.6,out=90,in=90] (m2.center) to (m1.center);
  \draw[thickstring] (m1.center) to (m1 |- b.north);
  \draw[thickstring] (m2.center) to (m2 |- bb.north);
  \draw[string] (in2.center) to (m3.center) to (m2 |- m3);
  \node (out) at (0,-0.75) {};
  \node (leftact) at (0,0) {};
  \draw[thickstring,arrow=0.8] (in1.center) to (in1 |- b.north);
\end{pic}
;\quad
\vcenter{\xymatrix@C1truecm@R0.5truecm{
   M\otimes M 
   \ar@<0.5ex>[r]^-{\beta}
  &
   R
   \ar@<0.5ex>[l]^-{\barbeta}
  \\
   M^*\otimes M^*
   \ar@<0.5ex>[r]^-{\barbeta}
  &
   R
   \ar@<0.5ex>[l]^-{\beta}
}}
\end{align}
We will study the properties of this isomorphism below in more detail.

  A (finite) algebra $A$ can be represented by a map into an endomorphism
ring $\End(M)\cong M\otimes M^*$. The algebra product is mapped
homomorphically onto the natural product of endomorphisms given by the
universal evaluation map, that is composition of endomorphisms. For our
purpose we use maps $h,y$, see~\eqref{three-multiplications}, such that
$y\circ h = 1_A$, that is we use faithful representations. In the light
of Wedderburn's theorem we may even assume, for simplicity, that
$A\cong \End(M)$, hence assuming $A$ is simple such that
$h\circ y = 1_{\End(M)}$. Now we can look `inside' the multiplication
in $A$ obtaining the left isomorphism in the next display.
\begin{align}\label{three-multiplications}
\begin{pic}
  \node (i1) at (0.125,1.5) {};
  \node (i2) at (1.375,1.5) {};
  \node[blackdot] (mul) at (0.75,0) {};
  \node (o) at (0.75,-1.5) {};
  \draw[string,out=270,in=180] (i1) to (mul.center);
  \draw[string,out=270,in=0] (i2) to (mul.center);
  \draw[string] (mul.center) to (o);
\end{pic}
\,\cong\,
\begin{pic}
  \node (i1) at (0.125,1.5) {};
  \node (in1) at (0,0.75) {};
  \draw[string] (i1) to (i1 |- in1);
  \node (in2) at (0.25,0.75) {};
  \node (in3) at (1.25,0.75) {};
  \node (i2) at (1.375,1.5) {};
  \draw[string] (i2) to (i2 |- in3);
  \node (in4) at (1.5,0.75) {};
  \node (m1) at (0,0.25) {};
  \node (m2) at (0.25,0.3) {};
  \node (m3) at (1.25,0.3) {};
  \node (m4) at (1.5,0.25) {};
  \node (d1) at (0.625,-0.5) {};
  \node (d2) at (0.875,-0.5) {};
  \node (out1) at (0.625,-0.75) {};
  \node (out2) at (0.875,-0.75) {};
  \node (o) at (0.75,-1.5) {};
  \draw[string] (o |- out1) to (o);
  \draw[thickstring,arrow=0.7] (in1) to (m1.center) .. controls +(0,-0.25) and +(0,.25) .. (d1.center) to (out1);
  \draw[thickstring,reverse arrow=0.7] (in4) to (m4.center) .. controls +(0,-0.25) and +(0,.25) .. (d2.center) to (out2);
  \draw[thickstring] (in2) to (m2.center);
  \draw[thickstring] (in3) to (m3.center);
  \draw[thickstring,arrow=0.55,out=270,in=270] (m3.center) to (m2.center);
  \node[coupon,fill=white] at (0.125,0.8) {h};
  \node[coupon,fill=white] at (1.375,0.8) {h};
  \node[coupon,fill=white] at (0.75,-0.8) {y};
\end{pic}
\,\cong\,
\begin{pic}
  \node (i1) at (0.125,1.5) {};
  \node (in1) at (0,0.75) {};
  \draw[string] (i1) to (i1 |- in1);
  \node (in2) at (0.25,0.75) {};
  \node (in3) at (1.25,0.75) {};
  \node (i2) at (1.375,1.5) {};
  \draw[string] (i2) to (i2 |- in3);
  \node (in4) at (1.5,0.75) {};
  \node (m1) at (0,0.25) {};
  \node (m2) at (0.25,0.3) {};
  \node (m3) at (1.25,0.3) {};
  \node (m4) at (1.5,0.25) {};
  \node (d1) at (0.625,-0.5) {};
  \node (d2) at (0.875,-0.5) {};
  \node (out1) at (0.625,-0.75) {};
  \node (out2) at (0.875,-0.75) {};
  \node (o) at (0.75,-1.5) {};
  \draw[string] (o |- out1) to (o);
  \draw[thickstring,arrow=0.7] (in1) to (m1.center) .. controls +(0,-0.25) and +(0,.25) .. (d1.center) to (out1);
  \draw[thickstring,arrow=0.7] (in4) to (m4.center) .. controls +(0,-0.25) and +(0,.25) .. (d2.center) to (out2);
  \draw[thickstring,->] (in2) to (m2.center);
  \draw[thickstring,->] (in3) to (m3.center);
  \draw[thickstring,out=270,in=270] (m3.center) to node[midway,blackdot,label=above:$\beta $]{} (m2.center);
%
  \node[coupon,fill=white] (c1) at (0.125,0.8) {$\mathsf{\sqcap} $};
  \node[coupon,fill=white] (c2) at (1.375,0.8) {$\mathsf{\sqcap} $};
  \node[coupon,fill=white] at (0.75,-0.8) {$\mathsf{\sqcup} $};
\end{pic}
;\quad
\vcenter{\xymatrix@C1truecm@R0.5truecm{
   A
   \ar@<0.5ex>[r]^-{\text{h}}
  &
   M\otimes M^*
   \ar@<0.5ex>[l]^-{y}
  \\
   A
   \ar@<0.5ex>[r]^-{\sqcap}
  &
   M\otimes M
   \ar@<0.5ex>[l]^-{\sqcup}
  \\ 
}}
\end{align}
  The second isomorphism is more subtle, as it involves the bilinear form.
This multiplication is called $\beta$-multiplication and operates on
$M\otimes M$. The choice of $\beta$ has to be compatible with the
morphisms $\sqcup,\sqcap$ in~\eqref{three-multiplications}. With 
$a,b\in A$ such that $a=\sum a_{ij} x_i\otimes y_j,$ $b=\sum b_{ij}
u_i\otimes v_j$ one obtains the multiplication
\begin{align}
  ab &= \sum a_{ij} x_i\otimes y_j\sum b_{ij} u_i\otimes v_j
      = \sum a_{ij}b_{lm} \beta( y_j,u_l) x_i\otimes v_m
\end{align}
  We remark here only, that the information flow in the endomorphism
ring situation is different (having upwards/back in time flow)
compared to the $\beta$-multiplication (related to a Frobenius algebra)
which has only downward information flow. This difference allows one
in quantum teleportation to \emph{choose} an entangled Bell state
(related to $\barbeta$) and make different Bell measurements (related to
$\beta$ up to unitary transformations), while the endomorphic situation
does not allow one this freedom. Also in linguistic models of meaning this
may have some implications, see subsection~\ref{information-flow}.

\subsection{Some notions from ring and module theory}
  All modules we are going to use are finitely generated projective (f.p.)
over a base field or ring. This is implied by the invertibility of the
Frobenius bilinear form (parastrophic matrix) hinging on a good duality
theory. This enables one to dualize algebra structures providing a
coalgebra structure, which fails in the general situation. Let $A$ be
an $R$-algebra with structure maps $m_A,\eta_A$. We denote by
$A^{op}$ the opposite algebra over the same $R$-module $A$, with the
opposite multiplication $m_A^{op} = m_A\circ \sigma$. It is useful
to introduce the \emph{enveloping algebra} $A^e=A\otimes A^{op}$,
\index{algebra!enveloping}\index{enveloping algebra}
which allows one to rewrite $A,A$-bimodules ${}_AM_A$ as $A^e$-left modules.

  A \emph{derivation} $D : A \rightarrow M$ is a linear operator from
\index{derivation}
the $A,A$-bimodule $A$ to the $A,A$-bimodule $M$, such that
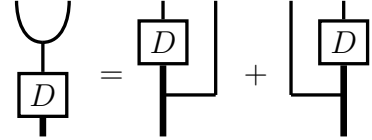
\begin{align}
 D(ab)&= D(a).b + a.D(b)
 &&&
\begin{pic}[scale=0.7]
  \node (in1) at (0,1.5) {};
  \node (in2) at (1,1.5) {};
  \node (mul) at (0.5,0.5) {};
  \node[coupon] (D) at (0.5,-0.5) {$D$};
  \node (out) at (0.5,-1.5){};
  \draw[string] (mul.center) to (D);
  \draw[thickstring] (D) to (out);
  \draw[string,out=270,in=180] (in1) to (mul.center);
  \draw[string,out=270,in=0] (in2) to (mul.center);
\end{pic}
\,\,=
\begin{pic}[scale=0.7]
  \node (in1) at (0,1.5) {};
  \node (in2) at (1,1.5) {};
  \node[coupon] (D) at (0,0.5) {$D$};
  \node (d1) at (1,-0.5) {};
  \node (out) at (0,-1.5) {};
  \draw[string] (in1) to (D);
  \draw[thickstring] (D) to (out);
  \draw[string] (in2) to  (d1.center);
  \draw[string] (d1.center) to (out |- d1);
\end{pic}
\,\,+
\begin{pic}[scale=0.7]
  \node (in1) at (0,1.5) {};
  \node (in2) at (1,1.5) {};
  \node (d1) at (0,-0.5) {};
  \node[coupon] (D) at (1,0.5) {$D$};
  \node (out) at (1,-1.5){};
  \draw[string] (in2) to (D);
  \draw[thickstring] (D) to (out);
  \draw[string] (in1) to (d1.center);
  \draw[string] (d1.center) to (out |- d1);
\end{pic}\end{align}
  where the module $M$ is represented by a bold line. The bold-unbold
`multiplication' like tangle is the right/left action of $A$ on $M$.
Let $\Der_R(A,M)$ be the $R$-module of derivations. A derivation $D_m$ is
called \emph{inner derivation} if there exists an $m\in M$ such that
\index{derivation!inner}
$D_m(a) = am-ma$. Now define the space of $A$-invariants of $M$ as
$M^A:=\{ m\in M \mid am=ma \}$. It is obvious that for all $m\in M^A$ the
inner derivation vanishes $D_m=0$. If $M=A$ as $R$-modules, the space of
invariants is just the kernel of the multiplication map $I(A)=\Ker(m_A)$.
One finds the following sequence to be exact
\begin{align}
  0 &\rightarrow M^A \rightarrow M \rightarrow \Der_R(A,M)
\end{align}
Using the isomorphisms between $A,A$-bimodules and $A^e$-left modules shows
that $M^A\cong \Hom_{A^e}(A,M)$ and $M\cong \Hom_{A^e}(A^e,M)$. It is also
easy to see that $m_A : A^e \rightarrow A$ is an epimorphism. Hence the
following sequence is exact
\begin{align}\label{exactseq}
  0 &\rightarrow I(A)=\Ker(m_a) \rightarrow A\otimes A^{op}
     \rightarrow A \rightarrow 0
\end{align}
A situation which is important in the Frobenius case is when this sequences
\index{split sequence}
is split. That is there exists a map $\delta : A\rightarrow I(A) :: a \mapsto
\delta(a)=a\otimes 1 - 1\otimes a$ whose image $I(A)$ in $A$ is an ideal
$A I(A)=I(A)=I(A) A$. Then $A^e=A\otimes A^{op}$ decomposes as a direct
sum $A^e=I(A)\oplus A$ and there is an idempotent pair $(e,1-e)$ projecting
\index{splitting idempotent}
onto the two spaces. This gives by standard algebra arguments some structure
results.
\mybenv{Lemma}
\begin{align}
\Hom_{A^e}(I(A),M) &\cong \Der_R(A,M)
\end{align}
\myeenv
Applying the functor $\Hom_{A^e}(-,A)$ to the exact sequence~\eqref{exactseq}
shows that $H\!H^1(A,$ $M) = \Ext^1_{A^e}(A,M) \cong \Der_R(A,M)/\InnDer_R(A,M)$,
where $H\!H^1$ is the first Hochschild cohomology group. As an aside, having a
\index{Hochschild cohomology}
Hopf algebra structure allows one to formalize several cohomology theories in a
uniform manner, see~\cite{sweedler:1968a}. The graphical calculus is not (very)
sensitive to the underlying ring structure, so we do not go deeper into ring
theory here. The main result we quote, establishes the existence of an spiting
idempotent for the class of finite projective algebras we are interested in.
\mybenv{Theorem}
Let $R$ be a commutative ring. For $R$-algebras $A$ the following statements
are equivalent:
\begin{itemize}
\item[(i)] $A$ is projective as a left $A^e$-module.
\item[(ii)] The exact sequence~\eqref{exactseq} for $A^e$-modules is split.
\item[(iii)] There exists a splitting idempotent element
     $e=\sum e_{(1)}\otimes e_{(2)} \in A\otimes A$ such that for all
     $a\in A,$ $ae=ea$ and $\sum e_{(1)}e_{(2)}=1$ holds.
\end{itemize}
\myeenv
An $R$-algebra $A$ is \emph{separable} iff $A/R$ is a separable ring extension.
\index{algebra!separable}
That is $m: A^e\rightarrow A$ is a split epimorphism of $A^e$-modules. By the
above theorem this is equivalent to say that $A$ is $A^e$-projective or that
there exists a splitting idempotent $e$ as in (iii). For further details
see~\cite[Sec. 5.2]{kadison:1999a}. The conditions in (iii) translates into the
following graphical statements ($\eta : I \rightarrow A$ is the unit map)
\begin{align}\label{splittingIdempotent}
\begin{pic}[scale=0.8]
  \node[blackdot] (top) at (0.5,1) {};
  \node (out1) at (0,0) {};
  \node (out2) at (1,0) {};
  \node[below] at (0.5,0) {$\Sigma e_{(1)}\otimes e_{(2)}$};
  \draw[string,out=180,in=90] (top) to (out1);
  \draw[string,out=0,in=90] (top) to (out2);
\end{pic}
:\,\,
\begin{pic}[scale=0.8]
  \node (in) at (0,1) {};
  \node (m1) at (0,0) {};
  \node (m2) at (1,0) {};
  \node (m3) at (2,0) {};
  \node[whitedot] (mul) at (0.5,-0.5) {};
  \node[blackdot] (top) at (1.5,0.5) {};
  \node (out1) at (0.5,-1) {};
  \node (out2) at (2,-1) {};
  \draw[string] (in) to (m1.center);
  \draw[string,out=270,in=180] (m1.center) to (mul);
  \draw[string] (mul) to (out1);
  \draw[string,out=180,in=90] (top) to (m2.center);
  \draw[string,out=270,in= 0] (m2.center) to (mul);
  \draw[string,out=0,in=90] (top) to (m3.center);
  \draw[string] (m3.center) to (out2);
\end{pic}
\,=
\begin{pic}[scale=0.8]
  \node (in) at (2,1) {};
  \node (m1) at (2,0) {};
  \node (m2) at (1,0) {};
  \node (m3) at (0,0) {};
  \node[whitedot] (mul) at (1.5,-0.5) {};
  \node[blackdot] (top) at (0.5,0.5) {};
  \node (out1) at (1.5,-1) {};
  \node (out2) at (0,-1) {};
  \draw[string] (in) to (m1.center);
  \draw[string,out=270,in=0] (m1.center) to (mul);
  \draw[string] (mul) to (out1);
  \draw[string,out=0,in=90] (top) to (m2.center);
  \draw[string,out=270,in= 180] (m2.center) to (mul);
  \draw[string,out=180,in=90] (top) to (m3.center);
  \draw[string] (m3.center) to (out2);
\end{pic}
\,=:
\begin{pic}[scale=0.8]
  \node (in) at (0.5,1) {};
  \node (m1) at (0,-0.5) {};
  \node (m2) at (1,-0.5) {};
  \node[whitedot] (mul) at (0.5,0) {};
  \node (out1) at (0,-1) {};
  \node (out2) at (1,-1) {};
  \draw[string] (in) to (mul);
  \draw[string,out=180,in=90] (mul) to (m1.center);
  \draw[string,out=0,in=90] (mul) to (m2.center);
  \draw[string] (m1.center) to (out1);
  \draw[string] (m2.center) to (out2);
\end{pic}
;
\begin{pic}[scale=0.8]
  \useasboundingbox (0,-1) rectangle (1,1); 
  \draw[line width=1.25pt] (0.5,0) circle (0.5cm);
  \node[blackdot] (top) at (0.5,0.5) {};
  \node[whitedot] (mul) at (0.5,-0.5) {};
  \node (out) at (0.5,-1) {};
  \draw[string] (mul) to (out);
\end{pic}
\,\,=
\begin{pic}[scale=0.8]
  \useasboundingbox (-0.3,-1) rectangle (0.3,1); 
  \node[whitedot] (eta) at (0,0.5) {$\eta $};
  \node (out) at (0,-1) {};
  \draw[string] (eta) to (out);
\end{pic}
\end{align}
  We defined the comultiplication map $\delta : A \rightarrow A\otimes A ::
a\mapsto ae=ea$. Coassociativity follows from the symmetric definition of
$\delta$ and from associativity of the product in $A$ and the `sliding'
of morphisms. Try it! For more information about splitting idempotents and
quadratic algebras see~\cite{hahn:1994a,caenepeel:militaru:zhu:2002a}, as
we want to avoid to discuss Azumaya and Taylor-Azumaya algebras. We will
use generators and bases so we quote two more standard results from algebra,
guaranteeing the existence of bases (generators).
\mybenv{Theorem}
  Any projective separable algebra $A$ over a commutative ring $R$ is finitely
generated. A separable algebra $A$ over a field $\Bbbk$ is semisimple.
\myeenv
Using this theorem, we find a finitely generated projective $R$-module $M$,
with generators $\{x_i\}_1^n$ and a dual module $M^*$ with dual basis
$\{f_i\}_1^n$ such that $A\cong M\otimes_R M^*$. The Frobenius homomorphism
will allow us to replace the dual module $A^*$ by $A$ and the dual basis by a
\emph{reciprocal basis} $\{y_i\}_1^n$.
\index{reciprocal basis}

\subsection{Frobenius functors}
\index{Frobenius!functors}
  Let $A,S$ be rings and let ${}_A\calM$, ${}_S\calM$ be the categories
of (f.p.) $A$-modules and $S$-modules. Let $i : S\rightarrow A$ be an 
injection, then any $A$-module can be turned into an $S$ module. This
defines the restriction functor $\sfR$ in the opposite direction
\begin{align}
  \sfR : {}_A\calM &\rightarrow {}_S\calM
  :: {}_AM\mapsto {}_SM :: sm = i(s)m 
\end{align}
  This is an instance of a change of base functor. Now $\sfR$ has a left 
adjoint $\sfT\dashv\sfR$ (induction functor) and a right adjoint 
$\sfR\dashv\sfH$ (coinduction functor) defined as follows.
\begin{align}
 \sfT : {}_S\calM & \rightarrow {}_A\calM ::
   \left\{\begin{array}{l}
            {}_SM \mapsto A\otimes_S {}_SM \\
            f \mapsto A\otimes_S f  
          \end{array}\right. \nn
 \sfH : {}_S\calM & \rightarrow {}_A\calM ::
   \left\{\begin{array}{l}
            {}_SM \mapsto  \Hom_S(A,{}_S\calM) \ni h  \\
            f \mapsto f\circ h  
          \end{array}\right.          
\end{align}
This means one has $\Hom_A(\sfT(M),N) \cong \Hom_S(M,\sfR(N))$ and 
$\Hom_S(\sfR(N),M) \cong \Hom(N,\sfH(M))$. The Frobenius property is
captured by the 
\mybenv{Definition}
A ring extension $A/S$ is a Frobenius extension if and only if
$\sfH$ and $\sfT$ are naturally equivalent as functors from
${}_S\calM \rightarrow {}_A\calM$.
\myeenv
We get a Frobenius
structure, see~\cite{kadison:1999a,caenepeel:militaru:zhu:2002a,khovanov:2010a},
\index{Frobenius!structure}
that is a triple $(\beta, \{x_i\},\{y_i\})$ where 
$\lambda\in \Hom_{S-S}(A,S)$ is the Frobenius homomorphism with
$\beta(a,b)=\lambda(ab)$ (inverse $\barbeta : S \rightarrow A\otimes A$),
and $\{x_i\},\{y_i\}$ are generators of $A$ fulfilling the
$\beta$-multiplication equations (see tangles~\eqref{betaInv}
and~\eqref{frobMove0}).
\begin{align}
  \sum_i x_i\beta(y_i,a) &= a,
  &&&
  \sum_i \beta(a,x_i)y_i = a
\end{align}
This is a generalization of the theorem~\ref{frobTheorem}, valid for
noncommutative rings, with many interesting applications, see for
example~\cite{khovanov:2004a}. Pairs of functors $(\sfT,\sfH)$
with $\sfT \dashv \sfR \dashv \sfH$ such that $\sfT\cong \sfH$ is an
isomorphism are called Frobenius pairs of Frobenius functors.

\subsection{Graphical characterization of Frobenius algebras}
Let ${}_R\calM$ be a tensor category of (f.p.) $R$-modules, and
$A\in {}_R\calM$. Let $A$ be an $R$-algebra with structure maps
($\mu_A,\eta_A$). If $A$ is Frobenius then further structure maps
exist, such as the Frobenius homomorphism $\lambda : A \rightarrow R$,
or equivalently an associative bilinear form
(see section~\ref{associativeBilinearForms})
$\beta = \lambda\circ\mu_A : A\otimes A \rightarrow R$. However,
it is graphically more effective to use the dual $\Lambda$
of the Frobenius homomorphism and the splitting idempotent element
$e$ to define a coalgebra structure $\delta : A\rightarrow A\otimes A$,
as in~\eqref{splittingIdempotent}. we define the coproduct as follows
\begin{align}\label{tang-FrobYanking}
\begin{pic}
  \node (in) at (0.5,0.75) {};
  \node (o1) at (0,-0.75) {};
  \node (o2) at (1,-0.75) {};
  \node[blackdot] (cmul) at (0.5,0) [label=north east:$\delta_A $] {};
  \draw[string] (in) to (cmul.center);
  \draw[string,out=180,in=90] (cmul.center) to (o1);
  \draw[string,out=  0,in=90] (cmul.center) to (o2);
\end{pic}
:=
\begin{pic}
  \node (in) at (1.5,0.75) {};
  \node[coupon,minimum width=0.75cm] (b) at (0.25,0.5) {$\barbeta $};
  \node (o1) at (0,-0.75) {};
  \node (o2) at (1,-0.75) {};
  \node[blackdot] (mul) at (1,0) [label=north:$m_A $] {};
  \node (m1) at (1.5,0.25) {};
  \node (m2) at (0.5,0) {};
  \draw[string] (in) to (m1.center);
  \draw[string,out=  0,in=270] (mul.center) to (m1.center);
  \draw[string,out=180,in=270] (mul.center) to (m2 |- b.south);
  \draw[string] (o1 |- b.south) to (o1);
  \draw[string] (mul.center) to (o2);
\end{pic}
=
\begin{pic}
  \node (in) at (0,0.75) {};
  \node[coupon,minimum width=0.75cm] (b) at (1.25,0.5) {$\barbeta $};
  \node (o1) at (1.5,-0.75) {};
  \node (o2) at (0.5,-0.75) {};
  \node[blackdot] (mul) at (0.5,0) [label=north:$m_A $] {};
  \node (m1) at (0,0.25) {};
  \node (m2) at (1,0) {};
  \draw[string] (in) to (m1.center);
  \draw[string,out=180,in=270] (mul.center) to (m1.center);
  \draw[string,out=  0,in=270] (mul.center) to (m2 |- b.south);
  \draw[string] (o1 |- b.south) to (o1);
  \draw[string] (mul.center) to (o2);
\end{pic}
;\qquad
\begin{pic}
  \node[coupon,minimum width=0.75cm] (b) at (0.25,0.5) {$\barbeta $};
  \node (o1) at (0,-0.75) {};
  \node[triangle] (t) at (0.5,-0.1) {$\lambda $};
  \draw[string] (o1 |- b.south) to (o1);
  \draw[string] (t |- b.south) to (t);
\end{pic}
\,=\,
\begin{pic}
  \node (i) at (0,-0.75) {};
  \node (o) at (0,0.75) {};
  \node[triangle,hflip] (t) at (0,0.1) {$\Lambda $};
  \draw[string] (i) to (t);
\end{pic}
\end{align}
It is easy to show graphically that $\Lambda$ is a unit for $m_A$ and
$\lambda$ is a counit for $\delta_A$. Hence we arrive at the
\mybenv{Definition}
\index{Frobenius algebra!graphical definition}
  A Frobenius algebra $A$ is a quintuple ($A,m_A,\lambda,\delta_A,\Lambda$)
such that $\lambda$ is a counit for $\delta_A$, $\Lambda$ is a unit for 
$m_A$ and the multiplication and comultiplication fulfill the
compatibility law
\begin{align}
\begin{pic}\label{tangleFrobeniusCompatibility}
  \node (i1) at (0,0.75) {};
  \node (i2) at (1.5,0.75) {};
  \node (o1) at (0.5,-0.75) {};
  \node (o2) at (2,-0.75) {};
  \node[blackdot] (m1) at (0.5,-0.25) {};
  \node[blackdot] (m2) at (1.5,0.25) {};
  \node (d1) at (0,0) {};
  \node (d2) at (1,0) {};
  \node (d3) at (2,0) {};
  \draw[string] (i1) to (d1.center);
  \draw[string] (m1.center) to (o1);
  \draw[string] (i2) to (m2.center);
  \draw[string] (d3.center) to (o2);
  \draw[string,out=270,in=180] (d1.center) to (m1.center);
  \draw[string,out=  0,in=270] (m1.center) to (d2.center);
  \draw[string,out= 90,in=180] (d2.center) to (m2.center);
  \draw[string,out=  0,in= 90] (m2.center) to (d3.center);
\end{pic}
=
\begin{pic}
  \node (i1) at (0,0.75) {};
  \node (i2) at (1,0.75) {};
  \node (o1) at (0,-0.75) {};
  \node (o2) at (1,-0.75) {};
  \node[blackdot] (m1) at (0.5,0.25) {};
  \node[blackdot] (m2) at (0.5,-0.25) {};
  \draw[string,out=270,in=180] (i1) to (m1.center);
  \draw[string,out=270,in=  0] (i2) to (m1.center);
  \draw[string] (m1.center) to (m2.center);
  \draw[string,out=180,in=90] (m2.center) to (o1);
  \draw[string,out=  0,in=90] (m2.center) to (o2);
\end{pic}
=
\begin{pic}
  \node (i1) at (0.5,0.75) {};
  \node (i2) at (2,0.75) {};
  \node (o1) at (0,-0.75) {};
  \node (o2) at (1.5,-0.75) {};
  \node[blackdot] (m1) at (0.5,0.25) {};
  \node[blackdot] (m2) at (1.5,-0.25) {};
  \node (d1) at (0,0) {};
  \node (d2) at (1,0) {};
  \node (d3) at (2,0) {};
  \draw[string] (i1) to (m1.center);
  \draw[string,out=180,in=90] (m1.center) to (d1.center);
  \draw[string,out=  0,in=90] (m1.center) to (d2.center);
  \draw[string,out=180,in=270] (m2.center) to (d2.center);
  \draw[string,out=  0,in=270] (m2.center) to (d3.center);
  \draw[string] (i2) to (d3.center);
  \draw[string] (o1) to (d1.center);
  \draw[string] (o2) to (m2.center);
\end{pic}  
\end{align}
together with the units and `yanking' rules given
in~\eqref{tang-FrobYanking} and their duals.
\myeenv
If we interpret multiplication as an $A$ action of the (left/right)
$A$-module $A$, and comultiplication similarly as (left/right) coaction,
then this compatibility relation reads `(left/right) actions and
coactions commute'.

\subsection{Algebraic characterizations of Frobenius algebras}
Frobenius algebras can be characterized in a number of ways, emphasizing
different aspects of this structure.
\index{Frobenius algebra!characterizations of}

\mybenv{Theorem}\cite{caenepeel:militaru:zhu:2002a}\label{frobTheorem}
Let $A$ be an $n$ dimensional $\Bbbk$-algebra, the following statements
are equivalent:
\begin{enumerate}
\item $A$ is Frobenius. 
\item There exists a \emph{Frobenius isomorphism}
   $\beta^r\in\Hom_\Bbbk({}_AA,{}_AA^*)$ for left $A$-modules in ${}_A\mathcal{M}$.
\item There exists a \emph{Frobenius isomorphism}
   $\beta^l\in\Hom_\Bbbk(A_A,A^*_A)$ for right $A$-modules in $\mathcal{M}_A$.
\item The left regular representation $\LR$ and right regular representation
   $\RR$ are equivalent.
\item There exists $a_i\in\Bbbk$ such that the parastrophic
   matrix~\eqref{parastrophic} is invertible.
\item There exists a nondegenerate associative bilinear form
   $\beta : A\times A\rightarrow \Bbbk$ (associativity :
   $\beta(ab,c)=\beta(a,bc)$ for all $a,b,c\in A$).
\item There exists a hyperplane in $A$ that does not contain a nonzero right
   ideal of $A$.
\item There exists a pair $(\lambda,\barbeta)$, called Frobenius pair, where
   $\lambda\in A^*$ is a \emph{Frobenius homomorphism} and $\barbeta :
   \Bbbk\rightarrow A\otimes A$ ($\barbeta=\Delta(1)=\sum\barbeta_{(1)}
   \otimes\barbeta_{(2)}$) such that for all $a\in A$
   \begin{align}
     a\barbeta &= \barbeta a
     &&&
     \sum \lambda(\barbeta_{(1)})\barbeta_{(2)}
        &=\sum \barbeta_{(1)}\lambda(\barbeta_{(2)})=1
   \end{align}
   where we used Heyneman-Sweedler notation for the comultiplication.
\end{enumerate}
\myeenv

\subsection{Some properties of Frobenius and Hopf algebras}
Frobenius algebras share some similarities with Hopf algebras, but also
exhibit different features. We will discuss the relation between Frobenius
and finite Hopf algebras in section~\ref{finHopfAlgebra}, moreover see
the chapters of~\cite{majid:2012a,vercruysse:2012a} in this book.

In quantum information theory it is appropriate to distinguish two extremal
cases.
\mybenv{Definition}
  A Frobenius algebra $A$ is called \emph{special} or
  \emph{trivially connected}, if the loop operator equals identity
  $l=m_A\circ\delta_A=1_A$ (see~\eqref{loop-fork}), or connected if $l$
  is invertible. A Frobenius alegbra is called \emph{totally disconnected}
  if the loop operator decomposes as $l=m_A\circ\delta_A= \Lambda\circ\lambda$.
\myeenv
\index{Frobenius alegbra!special}\index{Frobenius alegbra!disconnected}
The special or connected case (l.h.s. eqn.~\eqref{dissConnected}) shows
that by `yanking' one generates a tensor state which cannot be factored,
while the (r.h.s) equation shows that the totally disconnected case
produces a product state. 
\begin{align}\label{dissConnected}
\begin{pic}
  \node[coupon,minimum width=1.2cm] (b) at (0.5,1) {$\barbeta $};
  \node[blackdot] (c) at (0,0.25) {};
  \node[blackdot] (m) at (0,-0.5) {};
  \node (o1) at (0,-1) {};
  \node (o2) at (1,-1) {};
  \node (m1) at (-0.5,-0.125) {};
  \node (m2) at (0.5,-0.125) {};
  \draw[string] (c |-b.south) to (c.center);
  \draw[string] (o2 |- b.south) to (o2);
  \draw[string] (m.center) to (o1);
  \draw[string,out=180,in=90] (c.center) to (m1.center);
  \draw[string,out=  0,in=90] (c.center) to (m2.center);
  \draw[string,out=270,in=180] (m1.center) to (m.center);
  \draw[string,out=270,in=  0] (m2.center) to (m.center);
\end{pic}
=
\begin{pic}
  \node[coupon,minimum width=1.2cm] (b) at (0.5,1) {$\barbeta $};
  \node (o1) at (0,-1) {};
  \node (o2) at (1,-1) {};
  \draw[string] (o1 |- b.south) to (o1);
  \draw[string] (o2 |- b.south) to (o2);
\end{pic}
;\quad
\begin{pic}
  \node[coupon,minimum width=1.2cm] (b) at (0.5,1) {$\barbeta $};
  \node[blackdot] (c) at (0,0.25) {};
  \node[blackdot] (m) at (0,-0.5) {};
  \node (o1) at (0,-1) {};
  \node (o2) at (1,-1) {};
  \node (m1) at (-0.5,-0.125) {};
  \node (m2) at (0.5,-0.125) {};
  \draw[string] (c |-b.south) to (c.center);
  \draw[string] (o2 |- b.south) to (o2);
  \draw[string] (m.center) to (o1);
  \draw[string,out=180,in=90] (c.center) to (m1.center);
  \draw[string,out=  0,in=90] (c.center) to (m2.center);
  \draw[string,out=270,in=180] (m1.center) to (m.center);
  \draw[string,out=270,in=  0] (m2.center) to (m.center);
\end{pic}
=
\begin{pic}
  \node[coupon,minimum width=1.2cm] (b) at (0.5,1) {$\barbeta $};
  \node[triangle] (t1) at (0,0.55) {$\lambda $};
  \node[triangle,hflip] (t2) at (0,-0.75) {$\Lambda $};
  \node (o1) at (0,-1) {};
  \node (o2) at (1,-1) {};
  \draw[string] (t1 |- b.south) to (t1);
  \draw[string]  (t2) to (o1);
  \draw[string] (o2 |- b.south) to (o2);
\end{pic}
=
\begin{pic}
  \node (dummy) at (0,1.25) {};
  \node[triangle,hflip] (t1) at (0,0.25) {$\Lambda $};
  \node[triangle,hflip] (t2) at (1,0.25) {$\Lambda $};
  \node (o1) at (0,-1) {};
  \node (o2) at (1,-1) {};
  \draw[string] (t1) to (o1);
  \draw[string] (t2) to (o2);
\end{pic}
\end{align}
In \cite{coecke:pavlovic:vicary:2008a,coecke:paquette:pendrix:2008a} it
is shown that special Frobenius algebras are in one to one correspondence
with a choice of an orthonormal basis. Moreover one can characterize
classical structures and complementarity on q-bits using special Frobenius
algebras~\cite{coecke:duncan:2011a}.

\begin{align}\label{loop-fork}
\hskip-2ex
l &=
\begin{pic}
  \node (i) at (0.5,1) {};
  \node[blackdot] (c) at (0.5,0.5) {};
  \node[blackdot] (m) at (0.5,-0.5) {};
  \node (o) at (0.5,-1) {};
  \node (d1) at (0,0) {};
  \node (d2) at (1,0) {};
  \draw[string] (i) to (c.center) (m.center) to (o);
  \draw[string,out=180,in=90] (c.center) to (d1.center);
  \draw[string,out=270,in=180] (d1.center) to (m.center);
  \draw[string,out=0,in=90] (c.center) to (d2.center);
  \draw[string,out=270,in=0] (d2.center) to (m.center);
\end{pic}
;
\hskip1ex
f = \hskip-2ex
\begin{pic}
  \node (i1) at (0,1) {};
  \node (i2) at (1,1) {};
  \node[blackdot] (m) at (0.5,0.5) {};
  \node[blackdot] (c) at (0.5,-0.5) {};
  \node (o1) at (0,-1) {};
  \node (o2) at (1,-1) {};
  \draw[string] (c.center) to (m.center);
  \draw[string,out=270,in=180] (i1) to (m.center);
  \draw[string,out=270,in=0] (i2) to (m.center);
  \draw[string,out=180,in=90] (c.center) to (o1);
  \draw[string,out=0,in=90] (c.center) to (o2);
\end{pic}
;\hskip1ex
\begin{pic}
  \node (i) at (0,1) {};
  \node[coupon] (a) at (0,0) {$\alpha $};
  \node (o) at (0,-1) {};
  \draw[string] (i) to (a.north) (a.south) to (o);
\end{pic}
= 
\begin{pic}
  \node (i) at (0,1) {};
  \node[coupon,minimum width=0.6cm] (bb) at (0.6,0.6) {$\barbeta $};
  \node[coupon,minimum width=0.6cm] (b) at (0.6,-0.6) {$\beta $};
  \node (o) at (0,-1) {};
  \node (d1) at (0,0.4) {};
  \node (d2) at (0.4,0) {};
  \node (d3) at (0.8,0) {};
  \node (d4) at (0,-0.4) {};
  \draw[string] (i) to (d1.center) (d4.center) to (o);
  \draw[string] (bb.south -| d3) to (b.north -| d3);
  \draw[string,out=270,in=90] (d1.center) to (b.north -| d2);
  \draw[string,out=270,in=90] (bb.south -| d2) to (d4.center); 
\end{pic}  
;\hskip2.5ex
\begin{pic}
  \node (i1) at (0,1) {};
  \node at (0.5,0.75) {$\ldots $};
  \node (i2) at (1,1) {};
  \node[coupon,minimum width=1.2cm] (c) at (0.5,0) {foo};
  \node (o1) at (0,-1) {};
  \node at (0.5,-0.75) {$\ldots $};
  \node (o2) at (1,-1) {};
  \node at (0.5,-1.2) {\textrm{spider thm.}};
  \draw[string] (i1) to (i1 |- c.north) (i2) to (i2 |- c.north);
  \draw[string] (o1) to (o1 |- c.south) (o2) to (o2 |- c.south);
\end{pic}
=\hskip-1ex
\begin{pic}
  \node (i1) at (0,1.5) {};
  \node (i2) at (0.5,1.5) {};
  \node (i3) at (1.5,1.5) {};
  \node (o1) at (0,-1.5) {};
  \node (o2) at (0.5,-1.5) {};
  \node (o3) at (1.5,-1.5) {};
  \node[blackdot] (d1) at (0.25,1) {};
  \node (d2) at (1.5,1.25) {};
  \node[blackdot] (d3) at (0.625,0.8) [label=north east:$\ldots $] {};
  \node[blackdot] (d4) at (0.625,0.4) {};
  \node[blackdot] (d5) at (0.625,-0.1) {};
  \node[blackdot] (d6) at (0.625,-0.8) [label=south east:$\ldots $]{};
  \node[blackdot] (d7) at (0.25,-1) {};
  \node (d8) at (1.5,-1.25) {};
  \draw[string] (i3) to (d2.center) (d8.center) to (o3);
  \draw[string,out=270,in=180] (i1) to (d1.center);
  \draw[string,out=270,in=0] (i2) to (d1.center);
  \draw[string,out=270,in=180] (d1.center) to (d3.center);
  \draw[string,out=270,in=0] (d2.center) to (d3.center);
  \draw[string] (d3.center) to (d4.center);
  \draw[line width=1.25pt] (0.625,0.15) circle (0.25cm);
  \draw[dotted,line width=1.25pt] (d5.center) to (d6.center);
  \draw[string,out=0,in=90] (d6.center) to (d8.center);
  \draw[string,out=180,in=90] (d6.center) to (d7.center);
  \draw[string,out=180,in=90] (d7.center) to (o1);
  \draw[string,out=0,in=90] (d7.center) to (o2);
\end{pic}
\end{align}
  Let $l = m\circ\delta$ be a `loop'. One has the following normal
form (or spider) theorem
\index{Frobenius algebra!normal form theorem}
\index{spider theorem}
\mybenv{Theorem}
  Let $A$ be a symmetric Frobenius algebra, that is
  $m_A^{op} = m_A\circ\sigma=m_A$, then any tangle `foo' with arity
  $(n,m)$ ($n$ inputs, $m$-outputs) can be transformed using Frobenius
  moves and associativity to the normal form
  $\delta^{m-1} \circ l^r\circ m^{n-1}$ (with 
  $m^0=l=\delta^0=\Id$) for some non-negative integer $r$,
  see~\eqref{loop-fork}. If $A$ is special $l=1_A$.
\myeenv
This theorem is proven by recursion. The Frobenius
property~\eqref{tangleFrobeniusCompatibility} and as\-so\-cia\-ti\-vi\-ty
for $m_A$ and $\delta_a$ allows us to interchange the order of
multiplications and comultiplications.
In this way the inputs can be multiplied together, and the outputs
produced by comultiplications. In this process a certain number of `loops'
occur, which vanish if $A$ is special. Composing an $n,m$-tangle with
$k$ loops in its normal form with an $m,p$-tangle with $l$ loops
produces a $n,p$-tangle with $k+l+m$ loops. If $A$ is not symmetric
one needs to deal with the Nakayama automorphisms $\alpha$, see
section~\ref{sec:nakayama}. We note here that $\alpha$ can be
constucted using the bilinear forms $\beta,\barbeta$ and the trivial
symmetry $\sigma$ (switch) as in~\eqref{loop-fork} along the lines
we constructed the twist $\theta$ in~\eqref{Reidemeister2} from left/right
dualities and the braiding.

\mybenv{Theorem}
  A special symmetric Frobenius algebra $A$ is a bialgebra.
\myeenv
\begin{align}
&\hskip-2.5ex
\begin{pic}[scale=0.7]
  \node (i1) at (0,1.5) {};
  \node (i2) at (1,1.5) {};
  \node (o1) at (0,-1.5) {};
  \node (o2) at (1,-1.5) {};
  \node[blackdot] (c1) at (0,0.5) {};   
  \node[blackdot] (c2) at (1,0.5) {};
  \node[blackdot] (m1) at (0,-0.5) {};
  \node[blackdot] (m2) at (1,-0.5) {};
  \node (d1) at (-0.5,0) {};
  \node (d2) at (1.5,0) {};
  \draw[string] (i1) to (c1.center);
  \draw[string] (i2) to (c2.center);
  \draw[string] (o1) to (m1.center);
  \draw[string] (o2) to (m2.center);
  \draw[string,out=180,in=90] (c1.center) to (d1.center);
  \draw[string,out=270,in=180] (d1.center) to (m1.center);
  \draw[string,out=0,in=180] (c1.center) to (m2.center);
  \draw[string,out=180,in=0] (c2.center) to (m1.center);
  \draw[string,out=0,in=90] (c2.center) to (d2.center); 
  \draw[string,out=270,in=0] (d2.center) to (m2.center);
\end{pic}
=
\begin{pic}[scale=0.7]
  \node (i1) at (0,1.5) {};
  \node (i2) at (1,1.5) {};
  \node (o1) at (0.25,-1.5) {};
  \node (o2) at (1.25,-1.5) {};
  \node[blackdot] (c1) at (1,1) {};   
  \node[blackdot] (c2) at (1,0) {};
  \node[blackdot] (m1) at (1,0.5) {};
  \node[blackdot] (m2) at (0.25,-0.75) {};
  \node (d1) at (0,0.75) {};
  \node (d2) at (1.25,0.75) {};
  \node (d3) at (0,-0.4) {};
  \node (d4) at (1.25,-0.5) {};
  \draw[string] (i1) to (d1.center);
  \draw[string,out=270,in=180] (d1.center) to (m1.center);
  \draw[string] (i2) to (c1.center);
  \draw[string,out=0,in=90] (c1.center) to (d2.center);
  \draw[string,out=270,in=0] (d2.center) to (m1.center);
  \draw[string,out=180,in=0] (c1.center) to (m2.center);
  \draw[string] (m1.center) to (c2.center);
  \draw[string,out=180,in=270] (m2.center) to (d3.center);
  \draw[string,out=90,in=180] (d3.center) to (c2.center);
  \draw[string,out= 0,in=90] (c2.center) to (d4.center);
  \draw[string] (d4.center) to (o2);
  \draw[string] (m2.center) to (o1);   
\end{pic}
=
\begin{pic}[scale=0.7]
  \node (i1) at (0,1.5) {};
  \node (i2) at (1,1.5) {};
  \node (o1) at (0.5,-1.5) {};
  \node (o2) at (1.25,-1.5) {};
  \node[blackdot] (c1) at (1,1) {};   
  \node[blackdot] (c2) at (1,0) {};
  \node[blackdot] (m1) at (1,0.5) {};
  \node[blackdot] (m2) at (0.5,-0.5) {};  
  \node (d1) at (0,1) {};
  \node (d2) at (0,0) {};
  \node (d3) at (1.25,-0.5) {};
  \node (d4) at (1.25,0.75) {};
  \draw[string] (i1) to (d1.center);
  \draw[string,out=270,in=180] (d1.center) to (m1.center);
  \draw[string] (i2) to (c1.center);
  \draw[string,out=0,in=90] (c1.center) to (d4.center);
  \draw[string,out=270,in=0] (d4.center) to (m1.center);
  \draw[string] (m1.center) to (c2.center);
  \draw[string,out=180,in=90] (c1.center) to (d2.center);
  \draw[string,out=270,in=180] (d2.center) to (m2.center);
  \draw[string] (m2.center) to (o1);
  \draw[string,out=0,in=180] (m2.center) to (c2.center);
  \draw[string,out=0,in=90] (c2.center) to (d3.center);
  \draw[string] (d3.center) to (o2);
\end{pic}
=
\begin{pic}[scale=0.7]
  \node (i1) at (0,1.5) {};
  \node (i2) at (1,1.5) {};
  \node (o1) at (0,-1.5) {};
  \node (o2) at (1,-1.5) {};
  \node[blackdot] (c1) at (1,1) {};   
  \node[blackdot] (c2) at (0.5,-0.75) {};
  \node[blackdot] (m1) at (1,0.5) {};
  \node[blackdot] (m2) at (0.5,-0.25) {};  
  \node (d1) at (0,1) {};
  \node (d2) at (0,0.25) {};
  \node (d3) at (1.25,0.75) {};
  \draw[string] (i1) to (d1.center);
  \draw[string,out=270,in=180] (d1.center) to (m1.center);
  \draw[string] (i2) to (c1.center);
  \draw[string,out=0,in=90] (c1.center) to (d3.center);
  \draw[string,out=270,in=0] (d3.center) to (m1.center);
  \draw[string,out=180,in=90] (c1.center) to (d2.center);
  \draw[string,out=270,in=180] (d2.center) to (m2.center);
  \draw[string,out=270,in=0] (m1.center) to (m2.center);
  \draw[string](m2.center) to (c2.center);
  \draw[string,out=180,in=90] (c2.center) to (o1);
  \draw[string,out=0,in=90] (c2.center) to (o2);
\end{pic}
=
\begin{pic}[scale=0.7]
  \node (i1) at (0,1.5) {};
  \node (i2) at (1,1.5) {};
  \node (o1) at (0.75,-1.5) {};
  \node (o2) at (1.25,-1.5) {};
  \node[blackdot] (c1) at (1,1) {};   
  \node[blackdot] (c2) at (1,-1) {};
  \node[blackdot] (m1) at (0.375,0) {};
  \node[blackdot] (m2) at (1,-0.5) {};  
  \node (d1) at (0,0.75) {};
  \node (d2) at (0.75,0.75) {};
  \node (d3) at (1.25,0.75) {};
  \node (d4) at (0,0.25) {};
  \node (d5) at (0.75,0.25) {};
  \node (d6) at (1.25,0) {};
  \draw[string] (i1) to (d1.center);
  \draw[string,out=270,in=90] (d1.center) to (d5.center);
  \draw[string,out=270,in=90] (d2.center) to (d4.center);
  \draw[string] (i2) to (c1.center);
  \draw[string,out=180,in=90] (c1.center) to (d2.center);
  \draw[string,out=0,in=90] (c1.center) to (d3.center);
  \draw[string] (d3.center) to (d6.center);
  \draw[string,out=270,in=180] (d4.center) to (m1.center);
  \draw[string,out=270,in=0] (d5.center) to (m1.center);
  \draw[string,out=270,in=180] (m1.center) to (m2.center);
  \draw[string,out=270,in=0] (d6.center) to (m2.center);
  \draw[string] (m2.center) to (c2.center);
  \draw[string,out=180,in=90] (c2.center) to (o1);
  \draw[string,out=0,in=90] (c2.center) to (o2);
\end{pic}  
=
\begin{pic}[scale=0.7]
  \node (i1) at (0,1.5) {};
  \node (i2) at (1,1.5) {};
  \node (o1) at (0,-1.5) {};
  \node (o2) at (1,-1.5) {};
  \node[blackdot] (c1) at (0.5,0.5) {};   
  \node[blackdot] (c2) at (0.5,-1) {};
  \node[blackdot] (m1) at (0.5,1) {};
  \node[blackdot] (m2) at (0.5,-0.5) {};  
  \node (d1) at (0,0) {};
  \node (d2) at (1,0) {};
  \draw[string,out=270,in=180] (i1) to (m1.center);
  \draw[string,out=270,in=0] (i2) to (m1.center);
  \draw[string] (m1.center) to (c1.center);
  \draw[string,out=180,in=90] (c1.center) to (d1.center);
  \draw[string,out=270,in=180] (d1.center) to (m2.center);
  \draw[string,out=0,in=90] (c1.center) to (d2.center);
  \draw[string,out=270,in=0] (d2.center) to (m2.center);
  \draw[string] (m2.center) to (c2.center);
  \draw[string,out=180,in=90] (c2.center) to (o1);
  \draw[string,out=0,in=90] (c2.center) to (o2);
\end{pic}
=
\begin{pic}[scale=0.7]
  \node (i1) at (0,1.5) {};
  \node (i2) at (1,1.5) {};
  \node (o1) at (0,-1.5) {};
  \node (o2) at (1,-1.5) {};
  \node[blackdot] (c1) at (0.5,-0.5) {};   
  \node[blackdot] (m1) at (0.5,0.5) {};
  \draw[string,out=270,in=180] (i1) to (m1.center);
  \draw[string,out=270,in=0] (i2) to (m1.center);
  \draw[string] (m1.center) to (c1.center);
  \draw[string,out=180,in=90] (c1.center) to (o1);
  \draw[string,out=0,in=90] (c1.center) to (o2);
\end{pic}
\end{align}

Another interesting property, which can be used e.g. in singular value
decomposition~\cite{fauser:2004a}, is the following fact, true for any
convolution algebra $\Hom(A,A)$, hence for Hopf and Frobenius.
\mybenv{Theorem} The operators `loop' $l=m_A\circ \delta_A$ of arity
($1,1$) and `fork' $f=\delta_A\circ m_A$ of arity ($2,2$),
see~\eqref{loop-fork}, fulfill the same minimal polynomial, hence have
the same positive spectrum up to a null-space.
\myeenv
In the case of a special Frobenius algebra we see that $l=1_A$ and
hence $f$ is a projector onto a space isomorphic to $A$ in $A\otimes A$.

Kuperberg ladders~\cite{kuperberg:1991a,fauser:2002c} are the counterparts
in a Hopf algebra to the leftmost/rightmost tangles
in~\eqref{tangleFrobeniusCompatibility}. A Hopf algebra comes with an
antipode $\antip$~(\eqref{definitionHopf} and sec.~\ref{finHopfAlgebra}
below) causing the ladder tangle to be invertible (lhs
in~\eqref{kuperbergLadder}). These tangles play a role in invariant
theory of 3-manifolds.
\index{Kuperberg!ladder}
\begin{align}\label{kuperbergLadder}  
\begin{pic}[scale=0.7]
  \node (i1) at (0,1.5) {};
  \node (i2) at (1.5,1.5) {};
  \node (o1) at (1,-1.5) {};
  \node (o2) at (2.5,-1.5) {};kk
  \node[whitedot] (c1) at (1.5,1) {};
  \node[whitedot] (c2) at (2,0.5) {};
  \node[whitedot] (m1) at (0.5,-0.5) {};
  \node[whitedot] (m2) at (1,-1) {};
  \node (sb) at (1.5,0) {};
  \node (d1) at (0,0.5) {};
  \node (d2) at (2.5,0) {};
  \draw[string] (i1) to (d1.center);
  \draw[string,out=270,in=180] (d1.center) to (m1.center);
  \draw[string,out=0,in=180] (m1.center) to (c1.center);
  \draw[string] (i2) to (c1.center);
  \draw[string,out=0,in=90] (c1.center) to (c2.center);
  \draw[string,out=270,in=180] (m1.center) to (m2.center);
  \draw[string,out=0,in=270] (m2.center) to (sb.center);
  \draw[string,out=90,in=180] (sb.center) to (c2.center);
  \node[circle,draw,fill=white,inner sep=0mm] (S) at (1.5,-0.2) {$\antip $};
  \draw[string,out=0,in=90] (c2.center) to (d2.center);
  \draw[string] (d2.center) to (o2);
  \draw[string] (m2.center) to (o1);
  \node[whitedot] (c1) at (1.5,1) {};
  \node[whitedot] (c2) at (2,0.5) {};
  \node[whitedot] (m1) at (0.5,-0.5) {};
  \node[whitedot] (m2) at (1,-1) {};
\end{pic}
=
\begin{pic}[scale=0.7]
  \node (i1) at (0,1.5) {};
  \node (i2) at (1.5,1.5) {};
  \node (o1) at (0.5,-1.5) {};
  \node (o2) at (2,-1.5) {};
  \node[whitedot] (c1) at (1.5,1) {};
  \node[whitedot] (c2) at (1,0.5) {};
  \node[whitedot] (m1) at (1,-0.5) {};
  \node[whitedot] (m2) at (0.5,-1) {};
  \node (sb) at (1.5,0) {};
  \node (d1) at (0,-0.5) {};
  \node (d2) at (0.5,0) {};
  \node (d3) at (2,0.5) {};
  \draw[string] (i1) to (d1.center);
  \draw[string,out=270,in=180] (d1.center) to (m2.center);
  \draw[string] (m2.center) to (o1);
  \draw[string] (i2) to (c1.center);
  \draw[string,out=180,in=90] (c1.center) to (c2.center);
  \draw[string,out=180,in=90] (c2.center) to (d2.center);
  \draw[string,out=270,in=180] (d2.center) to (m1.center);
  \draw[string,out=270,in=0] (m1.center) to (m2.center);
  \draw[string,out=0,in=90] (c2.center) to (sb.center);
  \draw[string,out=270,in=0] (sb.center) to (m1.center);  
  \node[circle,draw,fill=white,inner sep=0mm] (S) at (1.5,0) {$\antip $};
  \draw[string,out=0,in=90] (c1.center) to (d3.center);
  \draw[string] (d3.center) to (o2);
  \node[whitedot] (c1) at (1.5,1) {};
  \node[whitedot] (c2) at (1,0.5) {};
  \node[whitedot] (m1) at (1,-0.5) {};
  \node[whitedot] (m2) at (0.5,-1) {};
\end{pic}
=
\begin{pic}[scale=0.7]
  \node (i1) at (0,1.5) {};
  \node (i2) at (0.5,1.5) {};
  \node (o1) at (0,-1.5) {};
  \node (o2) at (0.5,-1.5) {};
  \draw[string] (i1) to (o1);
  \draw[string] (i2) to (o2);
\end{pic}
;\quad
\begin{pic}
  \node (i1) at (0,1) {};
  \node (i2) at (1,1) {};
  \node[coupon] (h1) at (0,0.325) {$\sfH $};
  \node[coupon] (h2) at (1,0.325) {$\sfH $};
  \node[whitedot] (m) at (0.5,-0.25) [label=south east:$m_B$] {};
  \node (o) at (0.5,-1) {};
  \draw[string] (i1) to (h1.north);
  \draw[string,out=270,in=180] (h1.south) to (m.center);
  \draw[string] (i2) to (h2.north);
  \draw[string,out=270,in=0] (h2.south) to (m.center);
  \draw[string] (m.center) to (o);
  \node[whitedot,fill=gray] at (0.5,-0.25) {};
\end{pic}
=
\begin{pic}
  \node (i1) at (0,1) {};
  \node (i2) at (1,1) {};
  \node[coupon] (h) at (0.5,-0.325) {$\sfH $};
  \node[blackdot] (m) at (0.5,0.25) [label=north:$m_A$] {};
  \node (o) at (0.5,-1) {};
  \draw[string,out=270,in=180] (i1) to (m.center);
  \draw[string] (m.center) to (h.north);
  \draw[string,out=270,in=0] (i2) to (m.center);
  \draw[string] (m.center) to (h.north);
  \draw[string] (h.south) to (o);
\end{pic}
\end{align}
The rightmost equation in~\eqref{kuperbergLadder} shows a Frobenius
algebra homomorphism $\sfH : A \rightarrow B$ such that 
$m_B \circ (\sfH \otimes \sfH) = \sfH\circ m_A$. In red/green-calculus
the map in use is the Hadamard gate, which is invertible, allowing a
`color change', that is an change of algebra structure. As special
Frobenius algebras encode bases, this is essentially an entangling
operator changing the underlying classical structure. In the Hopf
algebraic case Sweedler developed a powerful cohomology
theory~\cite{sweedler:1968a} with a Hopf algebra action, which
provides a classification of such maps. Algebra homomorphisms fall
into the trivial cohomology class.

As a last example in this subsection we consider a module $A$
carrying a Hopf and a Frobenius algebra structure at the same time.
To make this situation well behaved we demand the following
\emph{distributive laws} (also called Laplace
property~\cite{rota:stein:1994a}) to hold as
compatibility relations. (White dots belong to the Hopf algebra,
black to Frobenius.)
\begin{align}
\begin{pic}[scale=0.8]
  \node (i1) at (0,1) {};
  \node (i2) at (1,1) {};
  \node (i3) at (1.5,1) {};
  \node (o) at (1,-1) {};
  \node (d) at (1.5,0.5) {};
  \node[whitedot] (w) at (0.5,0.5) {};
  \node[blackdot] (b) at (1,0) {};
  \draw[string,out=270,in=180] (i1) to (w.center);
  \draw[string,out=270,in=0] (i2) to (w.center);
  \draw[string,out=270,in=180] (w.center) to (b.center);
  \draw[string] (i3) to (d.center);
  \draw[string,out=270,in=0] (d.center) to (b.center);
  \draw[string] (b.center) to (o);
  \node[whitedot] at (0.5,0.5) {};
\end{pic}
=
\begin{pic}[scale=0.8]
  \node (i1) at (0,1.25) {};
  \node (i2) at (0.5,1.25) {};
  \node (i3) at (1.5,1.25) {};
  \node[whitedot] (w1) at (1.5,0.75) {};
  \node (d1) at (0,0.25) {};
  \node (d2) at (0.5,0.25) {};
  \node (d3) at (1,0.25) {};
  \node (d4) at (2,0.25) {};
  \node[blackdot] (b1) at (0.5,-0.25) {};
  \node[blackdot] (b2) at (1.5,-0.25) {};
  \node[whitedot] (w2) at (1,-0.75) {};
  \node (o) at (1,-1.25) {};
  \draw[string] (i1) to (d1.center);
  \draw[string,out=270,in=180] (d1.center) to (b1.center);
  \draw[string,out=0,in=270] (b1.center) to (d3.center);
  \draw[string,out=90,in=180] (d3.center) to (w1.center);
  \draw[string] (i3) to (w1.center);
  \draw[string,out=0,in=90] (w1.center) to (d4.center);
  \draw[string,out=270,in=0] (d4.center) to (b2.center);
  \draw[string,out=180,in=270] (b2.center) to (d2.center);
  \draw[string] (i2) to (d2.center);
  \draw[string,out=270,in=180] (b1.center) to (w2.center);
  \draw[string,out=0,in=270] (w2.center) to (b2.center);
  \draw[string] (w2.center) to (o);
  \node[whitedot] at (1.5,0.75) {};
  \node[whitedot] at (1,-0.75) {};
\end{pic}
;\quad
\begin{pic}[scale=0.8]
  \node (i1) at (0,1) {};
  \node (i2) at (0.5,1) {};
  \node (i3) at (1.5,1) {};
  \node (o) at (0.5,-1) {};
  \node (d) at (0,0.5) {};
  \node[whitedot] (w) at (1,0.5) {};
  \node[blackdot] (b) at (0.5,0) {};
  \draw[string,out=270,in=180] (i2) to (w.center);
  \draw[string,out=270,in=0] (i3) to (w.center);
  \draw[string,out=270,in=0] (w.center) to (b.center);
  \draw[string] (i1) to (d.center);
  \draw[string,out=270,in=180] (d.center) to (b.center);
  \draw[string] (b.center) to (o);
  \node[whitedot] at (1,0.5) {};
\end{pic}
=
\begin{pic}[scale=0.8]
  \node (i1) at (0.5,1.25) {};
  \node (i2) at (1.5,1.25) {};
  \node (i3) at (2,1.25) {};
  \node[whitedot] (w1) at (0.5,0.75) {};
  \node (d1) at (0,0.25) {};
  \node (d2) at (1,0.25) {};
  \node (d3) at (1.5,0.25) {};
  \node (d4) at (2,0.25) {};
  \node[blackdot] (b1) at (0.5,-0.25) {};
  \node[blackdot] (b2) at (1.5,-0.25) {};
  \node[whitedot] (w2) at (1,-0.75) {};
  \node (o) at (1,-1.25) {};
  \draw[string] (i1) to (w1.center);
  \draw[string,out=180,in=90] (w1.center) to (d1.center);
  \draw[string,out=270,in=180] (d1.center) to (b1.center);
  \draw[string,out=0,in=270] (b1.center) to (d3.center);
  \draw[string] (i3) to (d4.center);
  \draw[string,out=0,in=90] (w1.center) to (d2.center);
  \draw[string,out=270,in=0] (d4.center) to (b2.center);
  \draw[string,out=180,in=270] (b2.center) to (d2.center);
  \draw[string] (i2) to (d3.center);
  \draw[string,out=270,in=180] (b1.center) to (w2.center);
  \draw[string,out=0,in=270] (w2.center) to (b2.center);
  \draw[string] (w2.center) to (o);
  \node[whitedot] at (0.5,0.75) {};
  \node[whitedot] at (1,-0.75) {};
\end{pic}
\end{align}
In~\cite{fauser:2001b,brouder:fauser:frabetti:oeckl:2002a}, see
also~\cite{carroll:2005a}, it was demonstrated how this `Laplace Hopf
algebras' produce via twistings all multiplicative structures in a
perturbative quantum field theory. See
also~\cite{fauser:jarvis:2006a,fauser:2006a} for a fancy number
theoretic application.
In~\cite{fauser:jarvis:2003a,fauser:jarvis:king:wybourne:2005a,fauser:jarvis:king:2010a}
among others it was demonstrated how this structure underlies invariant
rings, providing powerful tools to shorten proofs and allowing to solve
otherwise difficult problems. As a mental picture, the Hopf algebra
operates on a tensor module as `concatenation', while the Frobenius
structure encapsulates the discrete permutation symmetries $S_k$ on
tensor terms $\sigma : \otimes^k A \rightarrow \otimes^k A$. 

\subsection{Associative regular bilinear forms and Frobenius homomorphisms}\label{associativeBilinearForms}
Frobenius algebras come with a second way, beside the closed structures, to
identify f.p. modules with their duals. This is essentially using the
parastrophic matrix from~\eqref{parastrophic}.
\mybenv{Definition}\cite{murray:2005a}
  Let $\Bbbk$ be a residual field. A \emph{regular associative bilinear form}
  is a $\Bbbk$-linear map $\beta\in\Bil(A,\Bbbk) : A\otimes_{\Bbbk} A
  \rightarrow \Bbbk$ such that
\index{bilinear form!regular, associative}
  \begin{align}
    \beta(ab,c) &= \beta(a,bc) &&&\textrm{associativity} \nn
    \forall a\in A\textrm{~with~}a\not=0,\,\,&\exists b\in A\,\textrm{such that}\,
    \beta(a,b)\not=0 &&&\textrm{non degenerate}
  \end{align}
  The linear form $\lambda:=\beta(-,1) = \beta(1,-)$ is called
  \emph{Frobenius homomorphism}. If $\lambda(ab)=\lambda(ba)$, that is
  \index{Frobenius!homomorphism}
  the Nakayama automorphisms $\alpha=\Id$ (see section~\ref{sec:nakayama}),
  it is called \emph{trace form} (see tangle in~\eqref{frobStructures}).
\index{trace form}
\myeenv
Usually bilinear forms are depicted as cup-tangles, to avoid confusion with
the evaluation maps from the closed structure we denote them as coupons.
\begin{align}\label{frobStructures}
\begin{pic}[scale=0.7]
  \node (in1) at (0,1) {};
  \node (in2) at (1,1) {};
  \node (in3) at (1.5,1) {};
  \node[blackdot] (mul) at (0.5,0) [label=south west:$m_A$] {};
  \node[coupon,minimum width=1cm] (b) at (1,-1) {$\beta $};
  \draw[string,out=270,in=180,->] (in1) to (mul);
  \draw[string,out=270,in=  0,->] (in2) to (mul);
  \draw[string] (mul) to (mul |- b.north);
  \draw[string,arrow=0.75] (in3) to (in3 |- b.north); 
\end{pic}
\,=
\begin{pic}[scale=0.7]
  \node (in1) at (0,1) {};
  \node (in2) at (0.5,1) {};
  \node (in3) at (1.5,1) {};
  \node[blackdot] (mul) at (1,0) [label=south east:$m_A$] {};
  \node[coupon,minimum width=1cm] (beta) at (0.5,-1) {$\beta $};
  \draw[string,out=270,in=180,->] (in2) to (mul);
  \draw[string,out=270,in=  0,->] (in3) to (mul);
  \draw[string] (mul) to (mul |- beta.north);
  \draw[string,arrow=0.75] (in1) to (in1 |- beta.north); 
\end{pic}
;\,
\begin{pic}[scale=0.7]
  \node (in1) at (0,1) {};
  \node[whitedot] (eta) at (1,0.325) {$\Lambda $};
  \node[coupon,minimum width=1cm] (beta) at (0.5,-1) {$\beta $};
  \draw[string,arrow=0.55] (eta) to (eta |- beta.north);
  \draw[string,arrow=0.55] (in1) to (in1 |- beta.north); 
  \node at (1.75,0) {$=$};
  \node (in2) at (2.5,1.25) {};
  \node[triangle] (lambda) at (2.5,-0.75) {$\lambda $};
  \draw[string,arrow=0.55] (in2) to (lambda);
  \node at (3,0) {$=$};
  \node (in3) at (4.75,1) {};
  \node[whitedot] (eta2) at (3.75,0.325) {$\Lambda $};
  \node[coupon,minimum width=1cm] (beta2) at (4.25,-1) {$\beta $};
  \draw[string,arrow=0.55] (eta2) to (eta2 |- beta2.north);
  \draw[string,arrow=0.55] (in3) to (in3 |- beta2.north); 
\end{pic}
\,;\,
\begin{pic}[scale=0.7]
  \node (in1) at (0,1) {};
  \node (in2) at (1,1) {};
  \node[coupon,minimum width=1cm] (beta) at (0.5,-1) {$\beta $};
  \draw[string,arrow=0.55] (in2) to (in2 |- beta.north);
  \draw[string,arrow=0.55] (in1) to (in1 |- beta.north); 
  \node at (1.5,0) {$=$};
  \node (in3) at (2,1) {};
  \node (in4) at (3,1) {};
  \node[blackdot] (mul) at (2.5,0.25) {};
  \node[triangle] (lambda) at (2.5,-0.5) {$\lambda $};
  \draw[string,out=270,in=180,->] (in3) to (mul);
  \draw[string,out=270,in=  0,->] (in4) to (mul);
  \draw[string,arrow=0.55] (mul) to (lambda);
\end{pic}
\end{align}
Two bilinear forms are related by a \emph{homothety},
\index{homothety!of bilinear forms}
$\beta\simeq \beta^\prime$ if there exists a unit $k\in\Bbbk^\times$ and an
automorphism $V\in \Aut_\Bbbk(A)$ such that
\begin{align}\label{bilinHomothety}
  \beta(a,b) &= k\beta^\prime(Va,Vb)
\end{align}
Two bilinear forms related by a homothety are not essentially different
and the set $\Bil$ can be partitioned into homothety equivalence classes.
The rightmost equation in~\eqref{frobStructures} shows further that the
Frobenius homomorphisms $\lambda$ and $\lambda^\prime = k\lambda V$ are
related if $V$ is an algebra homomorphism, as in the rhs
of~\eqref{kuperbergLadder}.

\subsection{Nakayama automorphism}\label{sec:nakayama}
We call a bilinear form symmetric if $\forall a,b\in A\,.\,
\beta(a,b)=\beta(b,a)$. Note that this does not imply that $A$ is symmetric,
in general $A\not=A^{op}$. The \emph{Nakayama automorphism}
\index{Nakayama automorphism}
$\alpha\in \Aut_{\Bbbk-alg}(A)$ measures the deviation from symmetry of
$\beta$.
\begin{align}
\beta(a,b) &= \beta(b,\alpha(a))
\end{align}
The Nakayama automorphism is unique up to inner automorphisms. First fix an
isomorphism $\phi : \Hom({}_AA,{}_AA^*)$. Any other such isomorphism is
given by first applying an automorphism to ${}_AA$ and then applying
$\phi$. This results in a transformation to a new Frobenius linear form
$\lambda^\prime =u\lambda : c \mapsto \lambda(ck)$ and a new bilinear form
$\beta^\prime(a,b) = \beta(a,bu)$, where $u\in A^\times$ is a unit in
the algebra $A$. The corresponding Nakayama automorphism transforms as
$\alpha^\prime = I_u\circ\alpha$, where $I_u(a) = uau^{-1}$ is an
inner automorphism.

The bilinear form $\beta$ is symmetric iff $\alpha=\Id$. Note that in a
braided setting the ordinary switch can be addressed as a virtual
crossing~\cite{kauffman:1999a} not encoding over/under information.
\begin{align}
\begin{pic}[scale=0.7]
  \node (in1) at (0,1) {};
  \node (in2) at (1,1) {};
  \node[coupon,minimum width=1cm] (beta) at (0.5,-1) {$\beta $};
  \draw[string,arrow=0.55] (in2.center) to (in2 |- beta.north);
  \draw[string,arrow=0.55] (in1.center) to (in1 |- beta.north); 
  \node at (1.5,0) {$=$};
  \node[below] at (1.5,-1.5) {symmetry};
  \node (in3) at (2,1) {};
  \node (in4) at (3,1) {};
  \node (m1) at (2,0) {};
  \node (m2) at (3,0) {};
  \node[coupon,minimum width=1cm] (beta2) at (2.5,-1) {$\beta $};
  \draw[string,arrow=0.8] (in3.center) to (m2.center) to (m2 |- beta2.north);
  \draw[string,arrow=0.8] (in4.center) to (m1.center) to (m1 |- beta2.north); 
\end{pic}
\,;\,
\begin{pic}[scale=0.7]
  \node (in1) at (0,1) {};
  \node (in2) at (1,1) {};
  \node[coupon,minimum width=1cm] (beta) at (0.5,-1) {$\beta $};
  \draw[string,arrow=0.55] (in2.center) to (in2 |- beta.north);
  \draw[string,arrow=0.55] (in1.center) to (in1 |- beta.north); 
  \node at (1.5,0) {$=$};
  \node[below] at (1.5,-1.5) {Nakayama automorphism};
  \node (in3) at (2,1) {};
  \node (in4) at (3,1) {};
  \node (m1) at (2,0.25) {};
  \node (m2) at (3,0.25) {};
  \node[coupon] (alpha) at (3,-0.1) {$\alpha $};
  \node[coupon,minimum width=1cm] (beta2) at (2.5,-1) {$\beta $};
  \draw[string,arrow=0.7] (in3.center) to (m2.center) to (alpha) to (alpha |- beta2.north);
  \draw[string,arrow=0.8] (in4.center) to (m1.center) to (m1 |- beta2.north); 
\end{pic}
;\,\,
\begin{pic}[scale=0.7]
  \node (in1) at (0,1) {};
  \node (in2) at (1,1) {};
  \node[coupon] (v) at (0,0.2) {$V$ };
  \node[coupon,minimum width=1cm] (beta) at (0.5,-1) {$\beta $};
  \draw[string,arrow=0.55] (in2.center) to (in2 |- beta.north);
  \draw[string,arrow=0.4,arrow=0.9] (in1.center) to (v) to (v |- beta.north); 
  \node at (1.5,0) {$=$};
  \node[below] at (1.5,-1.5) {transposition};
  \node (in3) at (2,1) {};
  \node (in4) at (3,1) {};
  \node[coupon] (v2) at (3,0.2) {$V^t$ };
  \node[coupon,minimum width=1cm] (beta2) at (2.5,-1) {$\beta $};
  \draw[string,arrow=0.4,arrow=0.9] (in4.center) to (v2) to (v2 |- beta2.north);
  \draw[string,arrow=0.55] (in3.center) to (in3 |- beta2.north); 
\end{pic}
\end{align}
The nondegenerate bilinear form $\beta : A\otimes A\rightarrow \Bbbk$
induces the notion of an  adjoint on $\End(A)$ via
$\beta(a,V^tb) := \beta(Va,b)$, which
we call \emph{transposition}. This transposition fulfills the usual properties
\index{bilinear form!transposition of endomap}
such as linearity, $(UV)^t=V^tU^t$, and $(U^{-1})^t = (U^t)^{-1}$. In the
presence of a nontrivial Nakayama automorphism transposition is in general
\emph{not} an involution.
\index{transposition!non-involutive}
\begin{align}
  \beta(a,V^{t^2} b) &= \beta(V^t a,b)
      = \beta(\alpha^{-1} b, V^t a) = \beta(V \alpha^{-1} b,a)
      = \beta(a, \alpha V\alpha^{-1} b),
\end{align}
hence we get $V^{t^2}=\alpha V\alpha^{-1}$. If the Frobenius algebra is
symmetric, then $\alpha=\Id$ and transposition is an involution.
\mybenv{Lemma}\cite{murray:2005a}
Bilinear forms $\beta$ and $\beta^\prime$ are homothetic iff $\exists
k\in\Bbbk^\times$ and $V\in \Aut_\Bbbk(A)$ such that $\rho_u = kV^tV\in
\Aut_\Bbbk(A)$.
\myeenv
The reader may compare $\rho_u$ to the definition of a positive operator
in quantum mechanics. 
In the same line of thought, we notice the following: Let $A$ be a
Frobenius $\Bbbk$-algebra with bilinear form $\beta$ and Nakayama
automorphism $\alpha$. As seen above, two homothetic forms
$\beta^\prime(a,b)=\beta(a,bu)$ are related by a unit $u\in A^\times$
with Nakayama's $\alpha^\prime = I_u\circ\alpha$. This shows that the order
of the Nakayama automorphism is independent of the choice of the form
in the homothety equivalence class. This can be used to define the
following norm function: Let $\alpha^n=1$ be a Nakayama automorphism
of finite order,
\index{norm!induced by Nakayama automorphism}
define the algebra norm $N_\alpha(a) := a\alpha(a)\ldots
\alpha^{n-1}(a)$. This norm can be interpreted as the evaluation of a
term $t^n$ at $a$ in the Ore ring of right twisted polynomials
$A[t,\alpha]$, 
see~\cite{lam:leroy:1988a,bueso:et.al:2003a,abramov:le:li:2005a}. 
If $\alpha^n=I_a$ then $\alpha(a)=a$ and one gets 
$(\alpha^\prime)^n = I_{N_\alpha(u)a}$.
In the involutive case, important for quadratic algebras and Clifford
algebras~\cite{hahn:1994a}, or more generally for $*$-algebras, one
defines such norms using `special elements' or directly using the
involution. In quantum theory one usually assumes an involutive
$*$-automorphism. 

The closed structures allow us to relate morphisms $f$ in $\Hom(A,B)$ to
dualized morphisms $f^*$ in $\Hom(B^*,A^*)$ etc. bending lines up or down.
That is using the topological move R0~\eqref{R0} or `yanking'.
Having the Frobenius bilinear form $\beta$ available, we have a second
possibility to bend lines, which this time produces maps in $\Hom(A,B^*)$
etc. We need first to define the inverse map $\barbeta : \Bbbk \rightarrow
A\otimes A$, which exists due to nondegeneracy of $\beta$. We should write
$\beta_A$ and $\barbeta_A$ for the components of the map on the category,
but to unclutter notation we drop these indices.
\begin{align}\label{betaInv}
\xymatrix@C1truecm{
   A
   \ar[r]^-{A\otimes \barbeta}
   \ar@/_2.5ex/[rr]_{A}
  &
   A\otimes A\otimes A
   \ar[r]^-{\beta\otimes A}
  &
   A
}
&&&
\begin{pic}[scale=0.7]
  \node (in) at (0,1) {};
  \node (out) at (2,-1) {};
  \node[coupon,minimum width=1cm] (beta) at (0.5,-0.75) {$\beta $};
  \node[coupon,minimum width=1cm] (barbeta) at (1.5,0.75) {$\barbeta $};
  \node (dummy) at (1,0) {};
  \draw[string,arrow=0.55] (in.center) to (in |- beta.north);
  \draw[string,reverse arrow=0.5] (dummy |- beta.north) to (dummy |- barbeta.south);
  \draw[string,arrow=0.55] (out |- barbeta.south) to (out.center);
\end{pic}
\,\,=\,
\begin{pic}[scale=0.7]
  \node (in) at (0,1) {};
  \node (out) at (0,-1) {};
  \draw[string,arrow=0.5] (in.center) to (out.center);
\end{pic}
\end{align}
$\barbeta$ is a right inverse. Using the Nakayama automorphism $\alpha$
we see that it is also a left inverse ($\baralpha=\alpha^{-1}$):
\begin{align}\label{frobMove0}
\begin{pic}[scale=0.7]
  \node (in) at (0,1) {};
  \node (out) at (0,-1) {};
  \draw[string,arrow=0.5] (in.center) to (out.center);
\end{pic}
\,=\,
\begin{pic}[scale=0.7]
  \node (in) at (0,1) {};
  \node (out) at (2,-1) {};
  \node[coupon,minimum width=1cm] (beta) at (0.5,-0.75) {$\beta $};
  \node[coupon,minimum width=1cm] (barbeta) at (1.5,0.75) {$\barbeta $};
  \node (dummy) at (1,0) {};
  \draw[string,arrow=0.55] (in.center) to (in |- beta.north);
  \draw[string,reverse arrow=0.5] (dummy |- beta.north) to (dummy |- barbeta.south);
  \draw[string,arrow=0.55] (out |- barbeta.south) to (out.center);
\end{pic}
\,=\,
\begin{pic}[scale=0.7]
  \node (in) at (0,1.5) {};
  \node (out) at (2,-1.5) {};
  \node[coupon] (a1) at (1,-0.4) {$\alpha $};
  \node[coupon] (a2) at (1, 0.3) {$\baralpha $};
  \node[coupon,minimum width=1cm] (beta) at (0.5,-1.25) {$\beta $};
  \node[coupon,minimum width=1cm] (barbeta) at (1.5,1.25) {$\barbeta $};
  \draw[string,arrow=0.55] (in.center) to (in |- beta.north);
  \draw[string] (a1 |- beta.north) to (a1) to (a2) to (a2 |- barbeta.south);
  \draw[string,arrow=0.55] (out |- barbeta.south) to (out.center);
\end{pic}
\,=\,
\begin{pic}[scale=0.7]
  \node (in) at (2,1.5) {};
  \node (out) at (0,-1.5) {};
  \node[coupon] (a1) at (2,-0.4) {$\alpha $};
  \node[coupon] (a2) at (0, 0.3) {$\baralpha $};
  \node (m1) at (1,0.5) {};
  \node (m2) at (1,-0.5) {};
  \node (m3) at (0,-0.5) {};
  \node (m4) at (2,0.5) {};
  \node[coupon,minimum width=1cm] (beta) at (1.5,-1.25) {$\beta $};
  \node[coupon,minimum width=1cm] (barbeta) at (0.5,1.25) {$\barbeta $};
  \draw[string,arrow=0.2] (in.center) to (m4.center) to (m2.center) to (m2 |- beta.north);
  \draw[string,reverse arrow=0.55] (a1 |- beta.north east) to (a1) to (a2) to (a2 |- barbeta.south);
  \draw[string,arrow=0.8] (m1 |- barbeta.south) to (m1.center) to (m3.center) to (out.center);
\end{pic}
\,=\,
\begin{pic}[scale=0.7]
  \node (in) at (2,1) {};
  \node (out) at (0,-1) {};
  \node[coupon,minimum width=1cm] (beta) at (1.5,-0.75) {$\beta $};
  \node[coupon,minimum width=1cm] (barbeta) at (0.5,0.75) {$\barbeta $};
  \node (dummy) at (1,0){};
  \draw[string,arrow=0.55] (in.center) to (in |- beta.north);
  \draw[string,reverse arrow=0.5] (dummy |- beta.north) to (dummy |- barbeta.south);
  \draw[string,arrow=0.55] (out |- barbeta.south) to (out.center);
\end{pic}
\end{align}
Here the crossings are `virtual' that is the switch $\sigma$. The moves
in~\eqref{betaInv} and~\eqref{frobMove0} should be compared with the
topological moves~\eqref{R0} for the the cup/cap tangles of the closed
structure. Here the left/right aspect is taken care of by the Nakayama
automorphism.

A main characterization of a Frobenius algebra in theorem~\ref{frobTheorem},
and one which generalizes, is given by the Frobenius isomorphism
\index{Frobenius!isomorphism} 
${}_AA \cong {}_AA^*$ of the left $A$-modules. The Frobenius bilinear form
and its dual allow us to construct such morphisms together with the closed
structures~\cite{fuchs:2006a,fuchs:stigner:2008a}. We define left/right
module maps $\beta^r\in\Hom(A,A^*)$, $\beta^l\in\Hom(A,{}^*A)$ and their
inverses. The left/right aspect refers to the closed structures involved.
\begin{align}\label{homBetaR}
\begin{pic}
  \node (in) at (0,1) {};
  \node[coupon] (b) at (0,0) {$\beta^r $};
  \node (out) at (0,-1) {};
  \draw[string,arrow=0.5] (in.center) node[anchor=south] {$A$} to (b);
  \draw[string,reverse arrow=0.5] (b) to (out.center) node[anchor=north] {$A^*$};
\end{pic}
\,\,:=
\begin{pic}
  \node (in2) at (0,1) {};
  \node[coupon,minimum width=0.6cm] (b2) at (0.25,0) {$\beta $};
  \node (d1) at (0.5,0.5) {};
  \node (d2) at (1,0.5){};
  \node (out2) at (1,-1) {};
  \draw[string,arrow=0.5] (in2.center) node[anchor=south] {$A$} to (in2 |- b2.north);
  \draw[string] (d1.center) to (d1 |- b2.north);
  \draw[reverse arrow=0.5,string,out=90,in=90] (d1.center) node[anchor=south] {$A$} to (d2.center) node[anchor=south] {$A^*$};
  \draw[reverse arrow=0.5,string] (d2.center) to (out2.center) node[anchor=north] {$A^*$};
\end{pic}
\,;\quad
\begin{pic}
  \node (in) at (0,1) {};
  \node[coupon] (b) at (0,0) {$\barbeta^r $};
  \node (out) at (0,-1) {};
  \draw[string,reverse arrow=0.5] (in.center) node[anchor=south] {$A^*$} to (b);
  \draw[string,arrow=0.5] (b) to (out.center) node[anchor=north] {$A$};
\end{pic}
\,\,:=
\begin{pic}
  \node (in2) at (0,1) {};
  \node[coupon,minimum width=0.6cm] (b2) at (0.75,0) {$\barbeta $};
  \node (d1) at (0,-0.5) {};
  \node (d2) at (0.5,-0.5){};
  \node (out2) at (1,-1) {};
  \draw[string,reverse arrow=0.5] (in2.center) node[anchor=south] {$A^*$} to (d1.center);
  \draw[string] (d2.center) to (d2 |- b2.south);
  \draw[reverse arrow=0.5,string,out=270,in=270] (d1.center) node[yshift=-2pt,anchor=north] {$A^*$} to (d2.center) node[yshift=-2pt,anchor=north] {$A$};
  \draw[arrow=0.5,string] (out2 |- b2.south) to (out2.center) node[anchor=north] {$A$};
\end{pic}
\,;\quad
\begin{pic}
  \node (in) at (0,1) {};
  \node[coupon] (b1) at (0,0.5) {$\barbeta^r $};
  \node[coupon] (b2) at (0,-0.5) {$\beta^r $};
  \node (out) at (0,-1) {};
  \draw[string,<-] (in.center) node[anchor=south] {$A^*$} to (b1);
  \draw[string,arrow=0.8] (b1) to (b2);
  \draw[string,reverse arrow=0.75] (b2) to (out.center) node[anchor=north] {$A^*$};
\end{pic}
\,=\,
\begin{pic}
  \node (in) at (0,1) {};
  \node (out) at (0,-1) {};
  \draw[string,reverse arrow=0.5] (in.center) node[anchor=south] {$A^*$} to (out.center) node[anchor=north] {$A^*$};
\end{pic}
\end{align}
and
\begin{align}\label{homBetaL}
\begin{pic}
  \node (in) at (0,1) {};
  \node[coupon] (b) at (0,0) {$\beta^l $};
  \node (out) at (0,-1) {};
  \draw[string,arrow=0.5] (in.center) node[anchor=south] {$A$} to (b);
  \draw[string,reverse arrow=0.5] (b) to (out.center) node[anchor=north] {${}^*\!A $};
\end{pic}
\,\,:=
\begin{pic}
  \node (in2) at (1,1) {};
  \node[coupon,minimum width=0.6cm] (b2) at (0.75,0) {$\beta $};
  \node (d2) at (0,0.5) {};
  \node (d1) at (0.5,0.5){};
  \node (out2) at (0,-1) {};
  \draw[string,arrow=0.5] (in2.center) node[anchor=south] {$A$} to (in2 |- b2.north);
  \draw[string] (d1.center) to (d1 |- b2.north);
  \draw[reverse arrow=0.5,string,out=90,in=90] (d1.center) node[anchor=south] {$A$} to (d2.center) node[anchor=south] {${}^*\!A $};
  \draw[reverse arrow=0.5,string] (d2.center) to (out2.center) node[anchor=north] {${}^*\!A $};
\end{pic}
\,;\quad
\begin{pic}
  \node (in) at (0,1) {};
  \node[coupon] (b) at (0,0) {$\barbeta^l $};
  \node (out) at (0,-1) {};
  \draw[string,reverse arrow=0.5] (in.center) node[anchor=south] {${}^*\!A $} to (b);
  \draw[string,arrow=0.5] (b) to (out.center) node[anchor=north] {$A$};
\end{pic}
\,\,:=
\begin{pic}
  \node (in2) at (1,1) {};
  \node[coupon,minimum width=0.6cm] (b2) at (0.25,0) {$\barbeta $};
  \node (d2) at (0.5,-0.5) {};
  \node (d1) at (1,-0.5){};
  \node (out2) at (0,-1) {};
  \draw[string,reverse arrow=0.5] (in2.center) node[anchor=south] {${}^*\!A $} to (d1.center);
  \draw[string] (d2.center) to (d2 |- b2.south);
  \draw[reverse arrow=0.5,string,out=270,in=270] (d1.center) node[yshift=-2pt,anchor=north] {${}^*\!A $} to (d2.center) node[yshift=-2pt,anchor=north] {$A$};
  \draw[arrow=0.5,string] (out2 |- b2.south) to (out2.center) node[anchor=north] {$A$};
\end{pic}
\,;\quad
\begin{pic}
  \node (in) at (0,1) {};
  \node[coupon] (b1) at (0,0.5) {$\barbeta^l $};
  \node[coupon] (b2) at (0,-0.5) {$\beta^l $};
  \node (out) at (0,-1) {};
  \draw[string,<-] (in.center) node[anchor=south] {${}^*\!A $} to (b1);
  \draw[string,arrow=0.8] (b1) to (b2);
  \draw[string,reverse arrow=0.75] (b2) to (out.center) node[anchor=north] {${}^*\!A $};
\end{pic}
\,=\,
\begin{pic}
  \node (in) at (0,1) {};
  \node (out) at (0,-1) {};
  \draw[string,reverse arrow=0.5] (in.center) node[anchor=south] {${}^*\!A $} to (out.center) node[anchor=north] {${}^*\!A $};
\end{pic}
\end{align}
The proof that the identity holds in~\eqref{homBetaR} and~\eqref{homBetaL}
requires \emph{both} the inverse Frobenius bilinear form~\eqref{frobMove0}
and the topological move for right/left duality~\eqref{R0}, and is left as
an easy exercise.

We close this section about the Nakayama automorphism and its implications
on `yanking' moves by showing that the left/right duality imposed by the
closed structures is related to the module homomorphisms $\beta^\bullet$
and $\barbeta^\bullet$ (with $\bullet=r$ or $l$, Frobenius isomorphisms)
\begin{align}\label{cupViaFrob}
\begin{pic}
  \node (i1) at (0,1) {};
  \node (i2) at (1,1) {};
  \node[coupon] (b1) at (0,0) {$\beta^r $};
  \node[coupon] (b2) at (1,0) {$\beta^l $};
  \node (d1) at (0,-0.5) {};
  \node (d2) at (1,-0.5) {};
  \draw[string,arrow=0.5] (i1.center) node[anchor=south] {$A$} to (b1);
  \draw[string,reverse arrow=0.5] (i2.center) node[anchor=south] {${}^*\!A $} to (b2);
  \draw[string] (b1) to (d1.center);
  \draw[string] (b2) to (d2.center);
  \draw[string,arrow=0.5,out=270,in=270] (d1.center) to (d2.center);
\end{pic}
\,=\,
\begin{pic}
  \node (i1) at (0,1) {};
  \node (i2) at (2.125,1) {};
  \node[coupon,minimum width=0.6cm] (b1) at (0.25,0) {$\beta $};
  \node[coupon,minimum width=0.6cm] (b2) at (1.5,0) {$\barbeta $};
  \draw[string,arrow=0.5] (i1.center) node[anchor=south] {$A$} to (i1 |- b1.north);
  \node (u1) at (0.5,0.5) {};
  \node (u2) at (0.875,0.5) {};
  \node (l1) at (0.875,-0.75) {};
  \node (l2) at (1.25,-0.75) {};
  \node (l3) at (1.75,-0.5) {};
  \node (l4) at (2.125,-0.5) {};
  \draw[string,arrow=0.5] (l4.center) to (i2.center) node[anchor=south] {${}^*\!A $};
  \draw[string] (u1.center) to (u1 |- b1.north);
  \draw[string] (l2 |- b2.south) to (l2.center);
  \draw[string,reverse arrow=0.5] (u2.center) to (l1.center);
  \draw[string,arrow=0.5,out=90,in=90] (u2.center) node[yshift=1mm,anchor=south] {$A^*$} to (u1.center) node[yshift=1mm,anchor=south] {$A$};
  \draw[string,arrow=0.5,out=270,in=270] (l2.center) node[yshift=-1mm,anchor=north] {$A$} to (l1.center) node[yshift=-1mm,anchor=north] {$A^* $};
  \draw[string,reverse arrow=0.5,out=270,in=270] (l4.center) to (l3.center);
  \draw[string] (l3.center) to (l3 |- b2.south);
\end{pic}
\,=\,
\begin{pic}
  \node (i1) at (0,1){};
  \node (i2) at (0.75,1) {}; 
  \node (d1) at (0,-0.25) {};
  \node (d2) at (0.75,-0.25) {};
  \draw[string,arrow=0.5] (i1.center) node[above] {$A$} to (d1.center);
  \draw[string,reverse arrow=0.5] (i2.center) node[above] {${}^*\!A $} to (d2.center);
  \draw[string,arrow=0.55,out=270,in=270] (d1.center) to (d2.center);
\end{pic}
\,=\,
\begin{pic}
  \node (i1) at (0,1){};
  \node (i2) at (0.75,1) {}; 
  \node (u1) at (0,0.0) {};
  \node (u2) at (0.75,0) {};
  \node (d1) at (0,-0.5) {};
  \node (d2) at (0.75,-0.5) {};
  \node[coupon] (theta) at (0,0.5) {$\theta $};
  \draw[string,arrow=0.5] (i1.center) node[above] {$A$} to (theta);
  \draw[string] (theta) to (u1.center);
  \draw[string,reverse arrow=0.55] (i2.center) node[above] {${}^*\!A $} to (u2.center);
  \draw[string,reverse arrow=0.5,out=270,in=270] (d1.center) node[yshift=-1mm,below] {${}^*\!A $} to (d2.center) node[yshift=-1mm,below] {$A$};
  \draw[string,out=270,in=90] (u1.center) to (d2.center);
  \draw[string,out=270,in=90,cross] (u2.center) to (d1.center);
\end{pic}
\end{align}
In~\cite{fuchs:2006a} it is further graphically shown, that the Nakayama
automorphism $\alpha = \barbeta^l\,\beta^r$ is actually an algebra
automorphism. Compare this form of $\alpha$ with the form for $\alpha$
given in~\eqref{loop-fork} which does not use the closed structure or
the braid.

\subsection{Finite Hopf algebras as Frobenius algebras}\label{finHopfAlgebra}
\index{Hopf algebra!finite}
  Hopf algebras will be discussed at length in other chapters of this
book~\cite{majid:2012a,vercruysse:2012a}. We will provide here only the
basic facts which relates them to Frobenius algebras. In
definition~\ref{definitionHopf} we saw that a Hopf algebra over $H$ is
at the same time an unital algebra $(H,\mu_H,\eta_H)$ and a counital
coalgebra $(H,\Delta_H,\epsilon_H)$ which are compatible by the Hopf
compatibility law which includes the switch map $\sigma$. The
multiplication $\mu_H$ extends to a multiplication $\mu_{H\otimes H}$.
\begin{align}
  \Delta_H \,\mu_H &= \mu_{H\otimes H} (\Delta_H\otimes \Delta_H)
    &&&
  \mu_{H\otimes H} := (\mu_H\otimes \mu_H)(\Id\otimes \sigma\otimes \Id)
\end{align}
Hence $\Delta_H$ is an algebra morphism $(H,\mu_H)\rightarrow (H\otimes H,
\mu_{H\otimes H})$.
A Hopf algebra unifies the concept of a (Lie-)group and a (Lie-)algebra
at the same time. An element $g\in H$ is called \emph{group like} if
\index{Hopf algebra!group like element}
\index{group like element}
$\Delta(g)=g\otimes g$, that is the diagonal action or `copying'. An element
$p\in H$ is called \emph{primitive} (or algebra like) if
\index{Hopf algebra!primitive element}
\index{primitive element}
$\Delta(p)=p\otimes 1+1\otimes p$. Primitive elements generate for example
the (universal enveloping) Lie algebra of a Lie group seen as Hopf algebra.
This analogy extends to the action of a Hopf algebra on a module $M$. Hence
we have a $H$-action $H\otimes M \rightarrow M$ and a $H$-coaction
$M\rightarrow H\otimes M$, which need to fulfill an analogue of the Hopf
compatibility law. In graphical terms this reads, using white nodes for
Hopf co/multiplications, as follows
\begin{align}\label{defHopf}
\begin{pic}[scale=0.7]
    \node (u) [whitedot] at (0,1) {};
    \node (d) [whitedot] at (0,-1) {};
    \node at ([yshift=-3mm]u.center) [right] {$\mu$};
    \node at ([yshift=3mm]d.center) [right] {$\Delta$};
    \node (ul) at (-1,2) {$H$};
    \node (ur) at (1,2) {$H$};
    \node (dl) at (-1,-2) {$H$};
    \node (dr) at (1,-2) {$H$};
    \draw[string] (u) to (d);
    \draw[string,out=90,in=180] (dl.north) to (d.west);
    \draw[string,out=90,in=0] (dr.north) to (d.east);
    \draw[string,out=270,in=180] (ul.south) to (u.west);
    \draw[string,out=270,in=0] (ur.south) to (u.east);
\end{pic}
=
\begin{pic}[scale=0.7]
    \node (uul) at (-1,2) {$H$};
    \node (uur) at (1,2) {$H$};
    \node (ul) [whitedot] at (-1,1) {};
    \node (ur) [whitedot] at (1,1) {};
    \node (dl) [whitedot] at (-1,-1) {};
    \node (dr) [whitedot] at (1,-1) {};
    \node (ddl) at (-1,-2) {$H$};
    \node (ddr) at (1,-2) {$H$};
    \node at ([yshift=3mm]ul.center) [right] {$\Delta$};
    \node at ([yshift=3mm]ur.center) [right] {$\Delta$};
    \node at ([yshift=-3mm]dl.center) [right] {$\mu$};
    \node at ([yshift=-3mm]dr.center) [right] {$\mu$};
    \node at (0,0) [right] {$\sigma$};
    \draw[string,out=180,in=180] (ul.west) to (dl.west);
    \draw[string,out=0,in=180] (ul.east) to (dr.west);
    \draw[string,out=0,in=0] (ur.east) to (dr.east);
    \draw[string,out=180,in=0] (ur.west) to (dl.east);
    \draw[string] (ul.north) to (uul.south);
    \draw[string] (ur.north) to (uur.south);
    \draw[string] (dl.south) to (ddl.north);
    \draw[string] (dr.south) to (ddr.north);
    \draw[mydashed] (-2,-1.5) rectangle (2,.3);
    \node at (2,-1.6) [right] {$\mu_{H\otimes H}$};
\end{pic}
\hskip-3ex;\,\, 
\begin{pic}[scale=0.7]
  \node (0) at (0,0) {$M$};
  \node (1) at (0,1) {};
  \node (2) at (-1,2) {$H$};
  \node (3) at (1,2) {$M$};
  \draw[thickstring] (0.north) to (1.center);
  \draw[thickstring, out=0, in=270] (1.center) to (3.south);
  \draw[string, out=180, in=270] (1.center) to (2.south);
\end{pic}
;
\begin{pic}[scale=0.7]
  \node (0) at (0,0) {$M$};
  \node (1) at (0,-1) {};
  \node (2) at (-1,-2) {$H$};
  \node (3) at (1,-2) {$M$};
  \draw[thickstring] (0.south) to (1.center);
  \draw[thickstring, out=0, in=90] (1.center) to (3.north);
  \draw[string, out=180, in=90] (1.center) to (2.north);
\end{pic}
\end{align}
The multiplication map $\mu_{H\otimes H}$ is depicted as the dashed box,
modules receive bold lines. The graphical description makes it clear, that
one can interchange the role of multiplication and comultiplication and we
see that $\mu_H$ is a morphism of coalgebras $(H\otimes
H,\Delta_{H\otimes H})$ and $(H,\Delta_H)$
\begin{align}
  \Delta_H \,\mu_H &= (\mu_H\otimes \mu_H)\Delta_{H\otimes H}
    &&&
  \Delta_{H\otimes H} := (\Id\otimes \sigma\otimes \Id)(\Delta_H\otimes \Delta_H)
\end{align}
just moving the dashed box in the tangle up. The compatibility law for a Hopf
action on a left $H$-co/module $M$ has a bold rightmost line in the lhs
of~\eqref{defHopf}.

The conditions under which a finite Hopf algebra over a commutative ring
$R$ is Frobenius were worked out by~\cite{larson:sweedler:1969a} for $R$ a
principal ideal domain, and by~\cite{pareigis:1971a} for general $R$ with
$\Pic[R]=0$ (The Abelian Picard group $\Pic[R]$ consists of the set of isoclasses
$[X]$ of linebundles $X$ over $R$ (i.e. $X^*\otimes X\cong R$) with $\otimes$
as multiplication $[X]+[X^\prime] = [X\otimes X^\prime]$ and 
$X\mapsto X^*=\Hom(X,R)$ as inverse). These results have consequences
for the existence and invertibility of the antipode, which in turn is
relevant for Hopf algebra cohomology and invariants of 3-manifolds.
So called integrals provide the main tool to prove these facts. Integrals
allow one to construct Frobenius homomorphisms and equivalently the
Frobenius bilinear form.
\index{Frobenius!homomorphism}

Let $\epsilon : H \rightarrow R$ be the augmentation map, which is an
algebra homomorphism $\epsilon(ab)=\epsilon(a)\epsilon(b)$.
\mybenv{Definition}
\index{Hopf algebra!left/right integral}
  A \emph{right (left) integral} is a $\mu_r \in H$ ($\mu_l \in H$)
  satisfying for all $h\in H$ the relation $h\mu_r = \epsilon(h)\mu_r$
  ($\mu_l h= \epsilon(h)\mu_l$). The space of all right (left) integrals
  is denoted as
  \begin{align}
    \int_H^r
      &:=\{\mu_r\in H\mid \forall h\in H.\, h\mu_r =\epsilon(h)\mu_r\}
    \\
  (\int_H^l
      &:=\{\mu_l\in H\mid \forall h\in H.\, \mu_l =\epsilon(h)\mu_l\})
    \nonumber
  \end{align}
\myeenv
Graphically integrals look like
\begin{align}\label{defIntegrals}
\begin{pic}
  \node (i1) at (0,0.5) {};
  \node (i2) at (1,0.5) {};
  \node (u1) at (0,0.25) {};
  \node (u2) at (1,0.25) {};
  \node[whitedot] (mul) at (0.5,-0.25) {};
  \node[whitedot] (e) at (0.5,-0.75) {$\epsilon $};
  \draw[string] (i1.center) to (u1.center);
  \draw[string] (i2.center) to (u2.center);
  \draw[string,out=270,in=180] (u1.center) to (mul);
  \draw[string,out=270,in=  0] (u2.center) to (mul);
  \draw[string] (mul) to (e);
  \node at (1.5,0) {$=$};
  \node (i3) at (2,0.5) {};
  \node (i4) at (2.75,0.5) {};
  \node[whitedot] (e3) at (2,-0.75) {$\epsilon $};
  \node[whitedot] (e4) at (2.75,-0.75) {$\epsilon $};
  \draw[string] (i3.center) to (e3);
  \draw[string] (i4.center) to (e4);
\end{pic}
;\quad
\begin{pic}
  \node (i1) at (0,0.75) {};
  \node[blackdot] (i2) at (1,0.25) {};
  \node (u1) at (0,0.25) {};
  \node (u2) at (1,0.25) {};
  \node[whitedot] (mul) at (0.5,-0.25) {};
  \node (out) at (0.5,-0.75) {};
  \draw[string] (i1.center) to (u1.center);
  \draw[string] (i2.center) node[above] {$\mu^r $} to (u2.center);
  \draw[string,out=270,in=180] (u1.center) to (mul);
  \draw[string,out=270,in=  0] (u2.center) to (mul);
  \draw[string] (mul) to (out.center);
  \node at (1.5,0) {$=$};
  \node (i3) at (2,0.75) {};
  \node (o3) at (2,-0.75) {};
  \node[whitedot] (i4) at (2,0.1) {$\epsilon $};
  \node[blackdot] (o4) at (2,-0.3) {};
  \draw[string] (i3.center) to (i4);
  \draw[string] (o4) node[anchor=west] {$\mu^r $} to (o3.center);
\end{pic}
;\quad
\begin{pic}
  \node (i1) at (1,0.75) {};
  \node[blackdot] (i2) at (0,0.25) {};
  \node (u1) at (1,0.25) {};
  \node (u2) at (0,0.25) {};
  \node[whitedot] (mul) at (0.5,-0.25) {};
  \node (out) at (0.5,-0.75) {};
  \draw[string] (i1.center) to (u1.center);
  \draw[string] (i2.center) node[above] {$\mu^l $} to (u2.center);
  \draw[string,out=270,in=   0] (u1.center) to (mul);
  \draw[string,out=270,in= 180] (u2.center) to (mul);
  \draw[string] (mul) to (out.center);
  \node at (1.5,0) {$=$};
  \node (i3) at (2,0.75) {};
  \node (o3) at (2,-0.75) {};
  \node[whitedot] (i4) at (2,0.1) {$\epsilon $};
  \node[blackdot] (o4) at (2,-0.3) {};
  \draw[string] (i3.center) to (i4);
  \draw[string] (o4) node[anchor=west] {$\mu^l $} to (o3.center);
\end{pic}
\end{align}
Let $H$ be a finitely generated projective module over a commutative ring with
$\Pic[R]=0$ underlying a Hopf algebra $H$. (The mild condition $\Pic[R]=0$
can be lifted by studying quasi Frobenius rings, which we do not pursue,
e.g.~\cite{nicholson:yousif:2003a}.)
The crucial property we need to construct a Frobenius homomorphism is the fact
hat $\int^r_H$ is a one dimensional module over $R$
\begin{align}
  \int_H^l\, H &\simeq H \simeq H\,\int_H^r
  &&&
  \textrm{with}\hskip2ex\int_H^l &\cong R \cong \int_H^r
\end{align}
As $\int_H^r\cong {}_RR$ as an $R$-module is one dimensional, it is an
invertible module. The following theorem, taken from~\cite{lorenz:2011a},
summarizes the work
of~\cite{larson:sweedler:1969a,pareigis:1971a,oberst:schneider:1973a}.
\mybenv{Theorem}
Let $H$ be a finite projectively generated Hopf algebra over the commutative
ring $R$, then:
\begin{itemize}
\item The antipode $\antip$ is bijective (has a linear inverse). This implies
  $\int^r_H = \antip(\int_H^l)$.
\item $H$ is a Frobenius $R$-algebra iff $\int_H^r\cong R$. This holds true
  if $\Pic[R]=0$. Moreover, if $H$ is Frobenius, then the dual Hopf algebra
  $H^*$ is Frobenius too.
\item Let $H$ be Frobenius. Then $H$ is symmetric iff
  \begin{itemize}
  \item[(i)] $H$ is unimodular (i.e. $\int_H^r=\int_H^l$), and
  \item[(ii)] $\antip^2$ is an inner automorphism of $H$.
  \end{itemize}
\end{itemize}
\myeenv
The existence of a one-dimensional $R$-module of right/left integrals
entails the construction of the Frobenius isomorphism
$\beta^\bullet : H \rightarrow H^*$ using the Hopf algebra structure
on $H$. The following tangle diagrams explain how to construct the
Frobenius homomorphism, and Frobenius isomorphisms, out of an unimodular
integral $\Lambda=\mu_l=\mu_r$ and the Hopf algebra structure:
\index{Frobenius!homomorphism!from an Hopf integral}
\begin{align}\label{FrobeniusHomfromIntegrals1}
\begin{pic}[scale=0.8]
  \node (i1) at (0,1) {$H^*$};
  \node (i2) at (0.5,1) {$H^*$};
  \node[coupon,minimum width=0.6cm] (bb) at (0.25,-0.5) {$\barbeta $};
  \draw[string] (i1) to (i1 |- bb.north);
  \draw[string] (i2) to (i2 |- bb.north); 
\end{pic}
&:=
\begin{pic}[scale=0.8]
  \node (i1) at (0,1) {$H^*$};
  \node (i2) at (0.5,1) {$H^*$};
  \node[triangle,hflip] (t) at (1.25,0.25) {$\Lambda $};
  \node[whitedot] (w) at (1.25,0) {};
  \node (d1) at (0,-0.5) {};
  \node (d2) at (0.5,-0.5) {};
  \node (d3) at (1,-0.5) {};
  \node (d4) at (1.5,-0.5) {};
  \draw[string] (t) to (w.center);
  \draw[string,out=0,in=90] (w.center) to (d4.center);
  \draw[string,out=270,in=270] (d4.center) to (d2.center);
  \draw[string,out=180,in=90] (w.center) to (d3.center);  
  \draw[string,out=270,in=270] (d3.center) to (d1.center);
  \draw[string] (i1) to (d1.center);
  \draw[string] (i2) to (d2.center);
  \node[whitedot] at (1.25,0) {};
\end{pic}
;\quad
\begin{pic}[scale=0.8]
  \node (i1) at (0,1) {$H^*$};
  \node (i2) at (0.5,1) {$H^*$};
  \node (i3) at (1,1) {$H$};
  \node[whitedot] (w) at (1,0.5) {};
  \node[coupon,minimum width=1.2cm] (bb) at (0.5,-0.5) {$\barbeta $};
  \node (d1) at (0.5,0.6) {};
  \node (d2) at (0.5,0.1) {};
  \node (d3) at (1,0.1) {};
  \node (d4) at (1.5,0.1) {};
  \draw[string] (i1) to (i1 |- bb.north);
  \draw[string] (i2) to (d1.center);
  \draw[string] (i3) to (w.center);
  \draw[string,out=270,in=90] (d1.center) to (d3.center);
  \draw[string] (d3.center) to (d3 |- bb.north);
  \draw[string,out=180,in=90] (w.center) to (d2.center);
  \draw[string,out=0,in=90] (w.center) to (d4.center);
  \draw[string] (d2.center) to (d2 -| i1);
  \draw[string] (d4.center) to (d3.center);
  \node[whitedot] at (1,0.5) {};
\end{pic}
=
\begin{pic}[scale=0.8]
  \node (i1) at (0,1) {$H^*$};
  \node (i2) at (0.5,1) {$H^*$};
  \node (i3) at (1,1) {$H$};
  \node[coupon,minimum width=0.7cm] (bb) at (0.25,-0.5) {$\barbeta $};
  \node[circle,draw,inner sep=0.2mm] (e) at (1,-0.5) {$\epsilon $};
  \draw[string] (i1) to (i1 |- bb.north);
  \draw[string] (i2) to (i2 |- bb.north);
  \draw[string] (i3) to (e); 
\end{pic}
;\quad
\begin{pic}[scale=0.8]
  \node at (0,1) {};
  \node[circle,draw,inner sep=0.2mm] (u) at (0,0.5) {$u$};
  \node (d1) at (0.5,0) {};
  \node (d2) at (1,0) {};
  \node (o) at (1,-1) {};
  \node[coupon,minimum width=0.6cm] (bb) at (0.25,-0.5) {$\barbeta $};
  \draw[string] (u) to (i1 |- bb.north);
  \draw[string] (d1.center) to (i2 |- bb.north);
  \draw[string,out=90,in=90] (d1.center) to (d2.center);
  \draw[string] (d2.center) to (o);
\end{pic}
:=
\begin{pic}[scale=0.8]
   \node[circle,draw,inner sep=0.2mm] (e) at (0,0.5) {$\epsilon^* $};
   \node at (0,1) {};
   \node (o) at (0,-1) {};
   \draw[string] (e) to (o);
\end{pic}
\\
\label{FrobeniusHomfromIntegrals2}
\begin{pic}
  \node (i) at (0,1) {};
  \node[circle,draw,inner sep=0.1mm] (s) at (0,0) {$\antip $};
  \node (o) at (0,-1) {};
  \draw[string] (i) to (s.north);
  \draw[string] (s.south) to (o);
\end{pic}
&:=
\begin{pic}
  \node (i1) at (1,1) {};
  \node[circle,draw,inner sep=0.1mm] (u) at (0.5,0.5) {$u$};
  \node[coupon,minimum width=0.6cm] (bb) at (0.25,-0.5) {$\barbeta $};
  \node (d1) at (1,0.25) {};
  \node (d2) at (0.5,0.25) {};
  \node (d3) at (-0.25,0) {};
  \node (d4) at (0,0) {};
  \node (o) at (-0.25,-1) {};
  \draw[string] (u) to (u |- bb.north);
  \draw[string] (i1) to (d1.center) to (d2.center);
  \draw[string] (d4.center) to (d4 |- bb.north);
  \draw[string] (d3.center) to (o);
  \draw[string,in=90,out=90] (d3.center) to (d4.center);
\end{pic}
;\hskip1truecm
\begin{pic}
  \node (i1) at (0,1) {};
  \node (i2) at (0.5,1) {};
  \node (d) at (0,0) {};
  \node (o) at (0.5,-1) {};
  \draw[thickstring,reverse arrow=0.25,reverse arrow=0.75] (i2) to (o);
  \draw[string] (i1) to (d.center) to (d -| i2);  
\end{pic}
:=
\begin{pic}
  \node (i1) at (0,1) {};
  \node (i2) at (0.5,1) {};
  \node (d1) at (0,0.5) {};
  \node (d2) at (0.5,0.5) {};
  \node[circle,draw,inner sep=0.1mm] (s) at (0.5,-0.25) {$\antip $};
  \node (d3) at (0,0) {};
  \node (d4) at (0.5,0) {};
  \node (d5) at (0.5,-0.75) {};
  \node (o) at (0,-1) {};
  \draw[string] (i1) to (d1.center);
  \draw[thickstring,reverse arrow=1] (i2) to (d2.center);
  \draw[string,out=270,in=90] (d1.center) to (d4.center);
  \draw[thickstring,out=270,in=90] (d2.center) to (d3.center);
  \draw[thickstring,reverse arrow=0.5] (d3.center) to (o);
  \draw[string] (d4.center) to (s.north);
  \draw[string] (s.south) to (d5.center);
  \draw[string] (d5.center) to (d5 -| o);
\end{pic}  
=
\begin{pic}
  \node (i1) at (0,1) {};
  \node (i2) at (0.5,1) {};
  \node (d1) at (0,0.5) {};
  \node (d2) at (0.5,0.5) {};
  \node[circle,draw,inner sep=0.1mm] (s) at (0.5,-0.25) {$\antip $};
  \node (d3) at (0,0) {};
  \node (d4) at (0.5,0) {};
  \node (d5) at (0.5,-0.6) {};
  \node (d6) at (0,-0.75) {};
  \node (d7) at (1,-0.75) {};
  \node (d8) at (1,0.1) {};
  \node (d9) at (1.5,0.1) {};
  \node (o) at (1.5,-1) {};
  \draw[string] (i1) to (d1.center);
  \draw[thickstring,reverse arrow=1] (i2) to (d2.center);
  \draw[string,out=270,in=90] (d1.center) to (d4.center);
  \draw[thickstring,out=270,in=90] (d2.center) to (d3.center);
  \draw[thickstring,reverse arrow=0.5] (d3.center) to (d6.center);
  \draw[string] (d4.center) to (s.north);
  \draw[string] (s.south) to (d5.center);
  \draw[string] (d5.center) to (d5 -| d7);
  \draw[thickstring,out=270,in=270] (d6.center) to (d7.center);
  \draw[thickstring,reverse arrow=0.6] (d7.center) to (d8.center);
  \draw[thickstring,out=90,in=90] (d8.center) to (d9.center);
  \draw[thickstring,reverse arrow=0.5] (d9.center) to (o);
\end{pic}
\end{align}
The first equation in~\eqref{FrobeniusHomfromIntegrals1} defines the
(dual) bilinear form from the integral $\Lambda$. The second equation
defines \emph{right orthogonality} of the form, see~\cite{larson:sweedler:1969a},
\index{bilinear form!orthogonality of an}
we need left orthogonality too. Orthogonality guarantees that the
Frobenius isomorphism respects the module structure. The third equation
defines the element $u$, which is needed to show the existence of the antipode
in~\eqref{FrobeniusHomfromIntegrals2}. Finally the right equation
in~\eqref{FrobeniusHomfromIntegrals2} defines the left action of $H$
on a dual module ${}_HM^*$ using the antipode and the right action on $M^*_H$,
which by duality~\eqref{leftrightaction} comes from the left action ${}_HM$,
establishing an Frobenius isomorphism.

\subsection{Information flow: Frobenius versus closed structures}\label{information-flow}
\index{information flow}\index{yanking}
A difference we hit over and over again is how the `yanking' is realized
in the closed and Frobenius situations. If we interpret the orientation
of tanlges as information flow, then a module propagates information in
time (downwards), while a dual module propagates information backwards
in time (upwards). A further difference is, that the closed structure
is characterized by a universal property, while the Frobenius structure
depends on a \emph{choice} of a bilinear form.

If we look at teleportation protocols modelled by `yanking', the ability
of Alice to chose between different Bell measurements indicates that she
is not dealing with a unique map, but with a Frobenius structure. Also
the creation of a shared entangled state is not unique, as every Bell
state does the trick. The need to communicate classical information to
Bob emerges from the need to communicate this choice of one of the
four Bell states in an \emph{a priory} mutually agreed on classical
basis, which can be described by a special Frobenius algebra. Hence the
`yanking' in teleportation should be thought of as a Frobenius algebra
related property.

A similar situation emerges in canonical quantization of fields. The negative
frequency parts are interpreted not as `particles' (field modes) travelling
backwards in time, but as `antiparticles' propagating homochronos. This
doubling of field modes resembles then the Frobenius type information flow
as discussed for teleportation. The $\beta$-multiplication is then given
via the reproducing kernel property of the field propagators
$\psi(x) =\int_Y g(x,y)\psi(y)$, 
but for a continuum. This calls for the extension of Frobenius structures
to a non-unital situation, as on an infinite dimensional space one
cannot have a unit, see~\cite{abramsky:heunen:2010a,coecke:heunen:2011a}.

With regard to the topic of this book, it is worth looking at the vector
space semantics of meanings of natural languages. There one has a set
of grammatical types of vector spaces. The meaning of words are
represented by vectors in differently typed vector spaces for nouns,
transitive verbs, etc. On top of these types one has a Lambeck 
pregroup, see~\cite{lambek:1999a,preller:sadrzadeh:2009a,sadrzadek:coecke:clark:2010a,preller:2012a}
and elsewhere in this volume, 
which is weakening the closed structure we are
using here. Using $a^l, a^r$ for left and right adjoints (duals), and
having an order on the `tensor monoid', one ends up, among others, with
relations of the type
$a^l a\leq 1 \leq aa^l$ and $aa^r \leq 1 \leq a^r a$. The later papers
in \emph{loc. cit.} show the striking similarity between this structure and
categorical quantum mechanics, especially teleportation. With respect to
our comment above, we think it should be investigated how this pregroup approach
relates to the Frobenius setup. First we note that any corpus of words is
finite, hence we can assume a good duality or a Frobenius structure on its
freely generated vector spaces of types. As the order of words in sentences is
crucial for their meaning, one has to deal with non-symmetric
Froebnius algebras, also to prevent left and right duality to coincide.
Hence one needs to take the Nakayama automorphism into account. As we
have seen in section~\ref{sec:nakayama} the braiding and the left/right
Frobenius isomorphisms can be transformed into one another~\eqref{cupViaFrob}.
The pictures in~\cite[page 150]{preller:sadrzadeh:2009a}
(and preface of this book)
describing sentences as \emph{John likes Mary} or \emph{John does not like Mary}
would then be replaced by tree like structures, where the multipications
are Frobenius pairings on the types, like $V\otimes W \rightarrow J$ for
a transitiv verb. The order structure may call for lax-Frobenius
structures. We close this speculation with the remark, that
in phylogenetic biomathematics one encounters very similar problems. Namely
to reconstruct ancestral relationships from present date gene sequences.
This is an analogue process to the linguistic setup, reconstructing a
tree which has as `words' the preset day gene expression and as 
`semantic meaning' the ancestral relation including branching times, see for 
example~\cite{jarvis:bashford:sumner:2005a,draisma:kuttler:2009a,jarvis:sumner:2011a,sumner:holland:jarvis:2011a}.

\section{A few pointers to further literature}
Frobenius algebras emerge in a large number of situations which we had no
opportunity to discuss here. We give a few hints where such developments
are found, and what problems they are addressing. Our pick on the literature
is subjective with an edge towards graphically minded work. 
\medskip

Starting with~\cite{abramsky:coecke:2004a} compact closed (dagger) categories
were used to analyze quantum theory, and especially quantum protocols.
Graphical methods have been utilized very heavily and led to `picturalism'
\index{picturalism}
in quantum theory~\cite{coecke:duncan:2011a}. The usage of \emph{several}
Frobenius algebra structures at the same time (red-green-/ black-white-/
rgb-calculi, which will be discussed elsewhere in this volume)
\index{graphical calculus!red-green}
\index{graphical calculus!black-white}
have captured, among other things, the concept of orthogonal
bases~\cite{coecke:pavlovic:vicary:2008a,coecke:paquette:pendrix:2008a},
\index{bases!orthogonal}
complementary bases, the dagger structure, and such things as connectedness
or disconnectedness of product states (e.g. entanglement for Bell, GHZ and
W states). 
\index{entanglement}
Specialness of Frobenius algebras is not found in invariant rings,
leading to non-trivial relations between reciprocal bases and dual bases, but still
Frobenius structures can be employed there~\cite{khovanov:1997a}. The product
of group characters can be understood along these lines, and the reader is
invited to stare at~\cite[2nd ed. pp305--309]{macdonald:1979a} to see how
it can be written in terms of the graphical calculus used here. A ring
extension $\mathbb{Q}[q,t]/\mathbb{Q}$ leads to Macdonald polynomials. Below
we will, however, give a pointer showing that it is not straight forward to
generalize this `Cartesian' setup (also needed for canonical quantization) to
manifolds and general coordinate invariance.
\medskip

A theme which is amenable to graphical treatment is the relation between
Frobenius structures of iterated ring extensions (or towers of algebras)
and the Jones polynomial, see~\cite{kadison:1999a,khovanov:1997a,mueger:2003a}.
\index{Jones polynomial}
If seen as ring extensions $A_n/A_{n-1}/\ldots /A_0$ the various
extensions provide separability idempotents, allowing the introduction
of the relevant Markov traces and finally of the Jones polynomial.
\medskip

A further theme is the relation of characters and certain trace modules
with Frobenius structures. A \emph{character} $\chi_V$ of a module $V$
in $A$-$\mathrm{mod}_R$ is given by
\begin{align}
  \chi_V(a) &= \mathrm{Tr}_{V/R}(a_V) \in R,\qquad (a\in A)
\end{align}
where $a_V\in\End(V)$ is given by $a_V(v)=av$. Hence characters form a subset
of trace forms
\begin{align}
  \chi_v\in \frac{A^*}{[A^*,A^*]}\subseteq A^* 
\end{align}
which vanish on $[A,A]$ (are constant on orbits). This allows one, for example,
to define the Higman trace $\tau=\tau_\beta$ for a Frobenius algebra $A$ as
\index{Higman trace}
\begin{align}
   \tau_\beta &: A \rightarrow A :: a \mapsto \sum_i x_i a y_i
\end{align}
independent on the choice of generators. $\tau$ is $\mathcal{Z}(A)$ linear,
with $\mathcal{Z}(A)$ the center of $A$. For a matrix algebra $A=M_n(R)$
one finds $\tau(a)=\mathrm{trace}(a)1_{n\times n}$, and for a group algebra
$A:=\mathbb{C}G$ one obtains the averaging or (up to normalization)
Reynolds operator $\tau(a)=\sum_{g\in G} gag^{-1}$.

Furthermore, one can define the Casimir operator, equivalent to the Higman
trace if $A$ is symmetric, as
\index{Casimir operator}
\begin{align}
  c &: A \rightarrow \mathcal{Z}(A) :: a \mapsto \sum_i y_i a x_i
\end{align}
which has deep connections to the Grothendieck groups, K-theory,
and restriction and induction functors, for example
see~\cite{lorenz:tokoly:2010a,lorenz:2011a} and references given there. It is
remarkable that the related diagrams, easily drawn, of the Higman trace, or the
Casimir operator emerge rather naturally in graphical calculations with
Frobenius algebras.
\medskip

As we mentioned already topological quantum field theory, we just remark that
Frobenius structures play a prominent role in TQFT,
see~\cite{atiyah:1989b,kock:2003a}. The general idea to apply a functor with
codomain a tensor monoid of vector spaces has proved to be very versatile.
Similar constructions can be found in the theory of vertex operators and rational
conformal quantum field theories, see for
example~\cite{fuchs:runkel:schweigert:2002a,fuchs:runkel:schweigert:2007a,barmeier.et.al:2010a},
which make heavy use of graphical calculi and provide further references.
\medskip

A theme related to classical physics, is that of Frobenius manifolds,
\index{Frobenius!manifolds}
see~\cite{hitchin:1997a}. Suppose $M$ is a manifold of dimension $n$, one can impose
the existence of the following sections
\begin{align}
  \theta \in C^\infty(T^*M) &&&&
  g \in C^\infty(S^2T^*M) &&&&
  c \in C^\infty(S^3T^*M)
\end{align}
on the (co)tangent space of $M$. Here $g$ is a metric on $M$, with covariant derivative
$\nabla$, and $\theta$ is a 1-form (related to the the Frobenius homomorphism)
and $c$ is a symmetric rank 3 tensor. Let $\{e_i\}$ be an orthogonal basis of $TM_x$,
with respect to the scalar product $(e_i,e_i)=\pm 1$, diagonalizing the left regular
representation $\LR_a$ for $a\in TM_x$. After rescaling, the $e_i$ are mutually
annihilating idempotents. Let $e=\sum_i e_i$ be the unit of this algebra, then we get
$\theta(v) = (e,v)$ and $(u,v)=c(e,u,v)$. Let $\mu_i=\theta(e_i)$ and let
$f_i\in TM^*_x$ be the dual basis of the $e_i$, then the structure maps can be written as
\begin{align}
  \theta =\sum \mu_i f_i &&&&
  g = \sum \mu_i f_i f_i &&&&
  c = \sum \mu_i f_i f_i f_i
\end{align}
The question which arises is, if this structure is compatible with the differential
structure, that is the covariant derivative, defined by the metric. This compatibility
leads to Chazy's \emph{non-linear} differential equation, see \emph{loc. cit.}. Chazy's
equation provides a `potential' for $\theta$, $g$ and $c$ resulting in the proper
definition of Frobenius manifolds.

If the metric is an Egoroff metric, having a sort of potential form,
$g= \sum_i \mu_i f_if_i = \sum_i\frac{\partial\phi}{\partial x_i}\textrm{d}x_i^2$, then
$\nabla c$ is symmetric, which implies $\textrm{d}_Ac=0$. This result relates orthogonal
coordinates (bases) to Frobenius structures, as we have seen in the quantum information
setting. However, on a manifold we need to be careful, as not all orthogonal frames are
allowed. This can be seen by the following example: The metric $\textrm{d}x^2+\textrm{d}y^2$
defines a Frobenius structure on $\mathbb{R}^2$, satisfying the condition above with
$\phi=x_1+x_2$. The same (in $\mathbb{R}^2\setminus\{0\}$) metric in polar coordinates
$\mathrm{d}r^2+r^2\mathrm{d}\phi^2$ does not carry a Frobenius structure with $\nabla c$
symmetric. This example shows clearly, that general coordinate transformations, hence
general covariance, are not compatible with (symmetric) Frobenius structures. It sheds
also some light on the usage of special Frobenius algebras in the semantics of quantum
protocols as mentioned above.
\medskip

Frobenius algebras can help to construct solutions of the Yang-Baxter
equation~\cite{baidar:fong:stolin:1997a} (which is related to our remark about
the Jones polynomial, and Kadison's work). In a boarder sense, one wants to
study entwined modules, that is, modules $A\otimes C$ where $A$ has an algebra
structure and $C$ has a coalgebra structure. It turns out, that Frobenius functors
play a crucial role in studying such entwined modules. This is developed at length
in~\cite{caenepeel:militaru:zhu:2002a}. It turns out, that this powerful
algebraic tools allows one to attack non-linear differential equations and provide
also solutions to the Yang-Baxter equation. A structure prominently used in
solid state and high energy physics problems often termed there `integrable models'.
\medskip

We have in this work always assumed that the Frobenius algebra structure is
\index{Frobenius algebra!non-associative}
associative. This can be relaxed to lax-Frobenius algebras or even general
non-associative Frobenius algebras. To my knowledge not much research has
been done in this direction, however see~\cite{oziewicz:wene:2011a}
where a first attempt is made to study non-associative Frobenius algebras,
also using graphical methods. Even in the seemingly trivial case of
complex numbers, which can be seen as a Frobenius algebra over the
reals~\cite{kock:2003a}, one encounters different Frobenius algebras
if associativity is dropped.

\medskip
\section*{Acknowledgement}
I would like to thank the organizers of the conference \emph{The categorical flow
of information in quantum physics and linguistics} for inviting me to this
interesting event. Special thanks to Chris Heunen for helping me with
the preparation of the \verb+TikZ+ graphics, and to Ronald C. King and
Peter D. Jarvis for a critical reading of a draft of this chapter.

%
%
%
\medskip
%

{\small
\bibliography{sql}
\bibliographystyle{agsm}
\def\topsep{0pt}
\def\parsep{0pt plus 5pt minus 1pt}
\def\itemsep{-0.5ex}
}
\end{document}